\documentclass[11pt]{amsart}
\usepackage{amsmath,amssymb,amsthm,mathtools,mathrsfs}
\usepackage[margin=1in]{geometry}
\usepackage{tikz}
\usetikzlibrary{arrows.meta}
\usepackage[hidelinks]{hyperref}
\mathtoolsset{showonlyrefs}
\allowdisplaybreaks[1]
\usepackage{cite}

\newtheorem{theorem}{Theorem}
\newtheorem{lemma}{Lemma}[section]
\theoremstyle{remark}
\newtheorem{remark}{Remark}[section]
\newcommand{\T}{\mathbb T}
\newcommand{\Z}{\mathbb Z}
\newcommand{\N}{\mathbb N}
\newcommand{\R}{\mathbb R}
\newcommand{\dd}{\,\mathrm d}

\usepackage{comment}

\newcommand{\leqc}{\lesssim}

\newcommand{\grad}{\nabla}

\newcommand{\dee}{\mathrm{d}}

\title[Osgood condition for anomalous dissipation]{Sharpness of the Osgood threshold for anomalous dissipation and non-uniqueness in passive scalar transport}
\author[K. Liss]{Kyle L. Liss}
\address[Kyle L. Liss]{Department of Mathematics, University of South Carolina, Columbia SC 29208, USA}
\email{kliss@mailbox.sc.edu}
\author[K. Luan]{Kunhui Luan}
\address[Kunhui Luan]{Department of Mathematics, University of South Carolina, Columbia SC 29208, USA}
\email{kunhui.luan@sc.edu}
\author[C. Tan]{Changhui Tan}
\address[Changhui Tan]{Department of Mathematics, University of South Carolina, Columbia SC 29208, USA}
\email{tan@math.sc.edu}
\date{}

\thanks{\textit{Acknowledgment.} 
CT is partially supported by NSF grant DMS-2238219.
}

\begin{document}

\begin{abstract}
We prove that for any concave, non-Osgood modulus of continuity $\omega$, there exists a divergence-free velocity field $u\in L^\infty([0,1];C^\omega(\T^2))$ and an initial data $f_0\in C^\infty(\T^2)$ such that the solutions of the associated advection--diffusion equation exhibit anomalous dissipation. This establishes the sharpness of the Osgood threshold for anomalous dissipation. It follows as a corollary of our main result that the Osgood condition is also sharp for the uniqueness of bounded solutions to the transport equation.
\end{abstract}

\maketitle

\section{Introduction}

We study the evolution of a passive scalar transported by an incompressible flow in the presence of weak molecular diffusion. Precisely, we consider the advection--diffusion equation
\begin{equation}\label{eq:ADE}
	\begin{cases}
		\partial_t f^\nu+u\cdot\nabla f^\nu=\nu\Delta f^\nu,\\
		f^\nu|_{t=0}=f_0\in L^\infty(\T^2),
	\end{cases}
\end{equation}
posed on $\T^2\cong[-\pi,\pi)^2$. Here, $u:[0,1]\times\T^2\to\R^2$ is a prescribed divergence-free vector field representing the velocity of a background fluid, $f^\nu:[0,1]\times\T^2\to\R$ is the passive scalar unknown, such as temperature or tracer concentration, and $0<\nu\ll1$ is the molecular diffusivity. The initial data $f_0:\T ^2 \to \R$ is bounded and fixed independently of $\nu$. Since \eqref{eq:ADE} preserves the spatial average of solutions, we may assume that $f_0$ is mean free. Then, $\|f^\nu(t)\|_{L^2}^2$ represents the scalar variance and may be viewed as an energy.

If $u\in L^\infty([0,1]\times\T^2)$, then for every $\nu>0$ there exists a unique solution
$$
f^\nu\in L^\infty([0,1]\times\T^2)\cap L^2([0,1];H^1(\T^2))
$$
satisfying the energy balance
\begin{equation}\label{eq:energyidentity}
	\frac12\|f^\nu(t)\|_{L^2}^2
	+\nu\int_0^t\|\grad f^\nu(s)\|_{L^2}^2\,\dee s
	=\frac12\|f_0\|_{L^2}^2.
\end{equation}
In particular, $t\mapsto\|f^\nu(t)\|_{L^2}$ is monotone decreasing and the integral
\begin{equation}\label{eq:dissipationdef}
	\mathcal E_\nu(t):=\nu\int_0^t\|\grad f^\nu(s)\|_{L^2}^2\,\dee s
\end{equation}
measures the dissipation of scalar energy over $[0,t]$. Formally, \eqref{eq:energyidentity} follows simply by pairing the equation with $f^\nu$, integrating by parts, and using that $u$ is divergence free.

Advection typically drives the formation of small spatial scales, which increases $\|\grad f^\nu(t)\|_{L^2}$ and consequently enhances the rate of scalar $L^2$ decay. When $u$ is sufficiently smooth, this effect is limited in the sense that $\mathcal E_\nu(t)\to0$ as $\nu\to0$ for any fixed $t>0$. Equivalently, energy is conserved in the vanishing diffusivity limit:
\begin{equation}\label{eq:noAD}
	\lim_{\nu\to0}\|f^\nu(t)\|_{L^2}=\|f_0\|_{L^2}.
\end{equation}
For example, if $u\in L^\infty([0,1];W^{1,\infty}(\T^2))$, then an energy estimate controls $\|\grad f^\nu\|_{L^\infty([0,1];L^2)}$ uniformly in $\nu$ when $f_0$ is smooth, and \eqref{eq:noAD} follows for bounded initial data by a standard approximation argument.

In this paper, we are interested instead in velocity fields that are continuous in space but rougher than Lipschitz. For such flows, advection can shift scalar Fourier mass to higher frequencies so rapidly that $\|\grad f^\nu(t)\|_{L^2}$ reaches a size comparable to $\nu^{-1/2}$ on a time scale independent of $\nu$. It then becomes possible that
\begin{equation}\label{eq:AD}
	\liminf_{\nu\to0}\mathcal E_\nu(t)>0
\end{equation}
for some $t>0$. This phenomenon is known as \textit{anomalous dissipation}.

Anomalous dissipation is a basic prediction of the Obukhov--Corrsin theory of passive scalar turbulence, which describes scalar statistics in a regime where the background turbulent flow exhibits H\"older-type scaling \cite{Obukhov49,Corrsin51,DSrY05,SS2000}. A dissipation anomaly is also central to theories of hydrodynamic turbulence. In three dimensions, the analogue of \eqref{eq:AD} is the prediction that kinetic energy dissipation remains nonvanishing as the viscosity tends to zero \cite{K41a,Frisch,EyinkSreenivasan2006}. This assumption is so fundamental to the theory that it is often called the ``zeroth law'' of turbulence.

Passive scalar transport provides a simpler setting than the Navier--Stokes equations for studying anomalous dissipation mathematically and has received considerable attention in recent years. The first deterministic example was obtained by Drivas, Elgindi, Iyer, and Jeong \cite{DEIJ22} using alternating shear flows that increase in frequency up to a finite singular time. For any fixed $\alpha\in(0,1)$, they construct a velocity field in $L^1([0,1];C^\alpha(\T^2))$ that leads to anomalous dissipation for a class of smooth initial data. A second example was given by Colombo, Crippa, and Sorella \cite{CCS}. They built upon finite-time mixing constructions \cite{Depauw2003,ACM,YaoZlatos2017} to obtain anomalous dissipation with uniform-in-time H\"older control on the velocity, while at the same time establishing the corresponding sharp scalar regularity predicted by the Obukhov--Corrsin theory.

 In both \cite{DEIJ22} and \cite{CCS}, the energy dissipation as $\nu \to 0$ is concentrated entirely at a terminal time where the velocity field loses smoothness. Armstrong and Vicol \cite{AV} introduced a new mechanism that produces continuous energy decay in the vanishing-diffusivity limit. The flow remains spatially rough throughout the evolution and sustains a cascade of scalar energy toward small scales. Their construction more closely reflects the physical picture of statistically stationary turbulence and provided the first example of anomalous dissipation by a single flow for \textit{all} sufficiently smooth initial data. Soon after, Burczak, Sz\'ekelyhidi, and Wu \cite{BurczakSzekelyhidiWu2023} combined their approach with convex integration to obtain anomalous dissipation for a velocity field that solves the Euler equation. In \cite{HessChildsRowan2025}, Hess-Childs and Rowan gave an example with total energy dissipation, in the sense that $\lim_{\nu \to 0}\|f^\nu(1)\|_{L^2}=0$ for all $L^2$ initial data. While the results above all use time-dependent flows, Johansson and Sorella \cite{JS2d} showed that anomalous dissipation can occur even for a two-dimensional autonomous flow with H\"{o}lder regularity. Scalar dissipation mechanisms have also been used to construct anomalous (kinetic energy) dissipation for the forced Navier--Stokes equations \cite{BDL,BCCDS,JS}. 

 In \cite{ElgindiLiss2024}, Elgindi and the first author of the present paper returned to the ideas of \cite{DEIJ22} and showed that alternating shear flows can be used to obtain anomalous dissipation for all smooth data, while also improving upon the spatial regularity of the flow.  In fact, the construction yields a single velocity field satisfying $u\in L^\infty([0,1];C^\alpha(\T^2))$ for \emph{every $\alpha<1$ simultaneously}. Since uniformly-in-time Lipschitz velocity fields cannot produce anomalous dissipation, this naturally raises the following question:

\medskip
\begin{center}
	\emph{What is the sharp spatial regularity threshold for anomalous dissipation to be possible?}
\end{center}
\medskip

\noindent This question is reminiscent of the classical problem of how far Lipschitz regularity of a vector field can be weakened while preserving uniqueness of the associated first-order ODE. The resemblance is not a coincidence, as anomalous dissipation is intimately related to properties of the underlying transport equation
\begin{equation}\label{eq:transportintro}
	\partial_t f+u\cdot\grad f=0.
\end{equation}
Indeed, if \eqref{eq:AD} holds, then sending $\nu \to 0$ and passing to a weak-$*$ limit of $\{f^\nu\}_{\nu>0}$ produces a bounded distributional solution of \eqref{eq:transportintro} that dissipates energy. It follows that anomalous dissipation is impossible whenever $u\in L^1([0,1];W^{1,1}(\T^2))$, since the DiPerna--Lions renormalization theory \cite{DL} guarantees that bounded transport solutions conserve their $L^2$ norm (see also \cite[Remark~4.3]{BonicattoCiampaCrippa2022}). More generally, the arguments in \cite[Theorem~3]{DEIJ22} and \cite[Theorem~1.3]{Rowan2024Kraichnan} rule out anomalous dissipation for velocity fields $u\in L^\infty([0,1]\times\T^2)$ whose associated transport equation has unique bounded solutions both forwards and backwards in time.

As in the ODE uniqueness problem, to describe sub-Lipschitz spatial regularity more precisely, it is natural to introduce the notion of modulus of continuity. Recall that a \emph{modulus of continuity} is a continuous, increasing function $\omega:[0,\pi]\to[0,\infty)$ satisfying $\omega(0)=0$. The corresponding spaces $C^\omega(\T^2)$ generalize H\"{o}lder spaces and are defined by the norm
$$
\|g\|_{C^\omega(\T^2)}
=\|g\|_{L^\infty(\T^2)}
+\sup_{0<d_{\T^2}(z,z')\le\pi}
\frac{|g(z)-g(z')|}{\omega(d_{\T^2}(z,z'))}.
$$
For a velocity field $u\in L^\infty([0,1];C^\omega(\T^2))$, it is well known that the classical \emph{Osgood condition} \cite{Osgood1898}
\begin{equation}\label{eq:Osgood}
	\int_0^\pi\frac{\dee r}{\omega(r)}=\infty
\end{equation}
is a sharp sufficient condition for uniqueness of trajectories for the ODE
\begin{equation}\label{eq:ODE}
    \frac{dX}{dt} = u(t,X(t)).
\end{equation}
In the divergence-free setting considered here, uniqueness of trajectories implies uniqueness of bounded solutions to \eqref{eq:transportintro} by Ambrosio's superposition principle \cite[Theorem~3.4]{AmbrosioCrippa2014}. Therefore, the Osgood condition rules out anomalous dissipation and gives an obstruction beyond Lipschitz regularity, including, for example, log-Lipschitz velocity fields.

The velocity fields from \cite{ElgindiLiss2024} could be designed to satisfy $u \in L^\infty([0,1];C^\omega(\T^2))$ with
\begin{equation}\label{eq:oldregularity}
\omega(r) = r|\log r|^p
\end{equation}
for any fixed $p>3$. These moduli are necessarily \textit{non-Osgood}, meaning that 
\begin{equation}\label{eq:nonOsgood}
	\int_0^\pi\frac{\dee r}{\omega(r)}<\infty.
\end{equation}
Thus, even though the construction reached every H\"older exponent below one, a gap remained between the regularity achieved and the Osgood obstruction.

It is natural to conjecture that the Osgood condition gives the sharp regularity threshold for anomalous dissipation. This was suggested in \cite[Section 3.1]{ElgindiLiss2024}, where it was asked whether given any non-Osgood modulus one can construct $u\in L^\infty([0,1];C^\omega(\T^2))$ that exhibits anomalous dissipation. Our main result answers this question affirmatively.

\begin{theorem}\label{thm:main}
Let $\omega$ be a concave and non-Osgood modulus of continuity. There exists a divergence-free velocity
\[
u\in L^\infty([0,1];C^\omega(\T^2))
\cap L^\infty_{\mathrm{loc}}([0,1);W^{1,\infty}(\T^2))
\]
and an integer $m\ge1$ such that the solutions $\{f^\nu\}_{\nu>0}$ of \eqref{eq:ADE} with $f_0(x,y)=\cos(mx)$ satisfy
\[
\liminf_{\nu\to0}\nu\int_0^1\|\nabla f^\nu(t)\|_{L^2}^2\dd t>0.
\]
\end{theorem}

\begin{remark}
    There is nothing particularly special about the choice $f_0(x,y) = \cos(mx)$. The same proof extends to a large class of initial data with Fourier support in a region near one of the coordinate axes of $\Z^2$. The resulting velocity fields, however, depend on the initial data. Since our primary interest was in the sharp regularity threshold, we did not pursue extensions to data-independent constructions here.
\end{remark}

\begin{remark}
An additional assumption on the modulus, such as concavity, is necessary. Indeed, there exists a non-Osgood modulus of continuity $\omega$ for which any velocity field $v\in L^\infty([0,1];C^\omega(\T^2))$ is actually log-Lipschitz uniformly in time, and hence cannot produce anomalous dissipation. Concavity rules out such pathological examples. The conclusion of Theorem~\ref{thm:main} holds also for sub-additive moduli, since they always have an equivalent concave majorant.
\end{remark}

Theorem~\ref{thm:main} also settles the corresponding sharp regularity problem for the uniqueness of bounded transport solutions with divergence-free drift. It has been known for some time that the Osgood condition is a sharp uniqueness criterion for nonnegative measure-valued solutions of the continuity equation
\begin{equation}\label{eq:continuity}
\partial_t\mu+\operatorname{div}(u\mu)=0.
\end{equation}
This is simply because PDE uniqueness is equivalent to ODE uniqueness in the setting of nonnegative measures. Indeed, non-unique trajectories immediately produce distinct solutions of \eqref{eq:continuity} by placing a Dirac mass along each trajectory, while the superposition principle \cite{Ambrosio2004,AmbrosioCrippa2014} supplies the other direction. Interestingly, this equivalence fails for general signed measures \cite{BonicattoGusev2019}. However, the stronger assumption of an Osgood modulus guarantees uniqueness even within this larger class of solutions \cite{AmbrosioBernard2008}.

One can then ask whether the Osgood condition remains necessary for uniqueness when the measures are required to be absolutely continuous. Colombo and Kumar \cite{ColomboKumar2026} recently showed that this is indeed the case, even under the restriction that the drift be divergence free. Given any concave and non-Osgood modulus $\omega$, they produced a divergence-free vector field $u \in C([0,1];C^\omega)$ such that \eqref{eq:transportintro} admits two distinct nonnegative solutions in $L^\infty([0,1];L^1)$, both starting from the characteristic function of a cube. Their construction involves concentration of the transported density, which does not remain uniformly bounded. Our result shows that the Osgood threshold remains sharp even when the solutions are required to be bounded. Precisely, combining Theorem~\ref{thm:main} with the time-reflection argument used in \cite[Theorem~3]{DEIJ22} gives the following.
\begin{theorem}\label{thm:nonuniqueness}
Let $\omega$ be a concave and non-Osgood modulus of continuity. There exists a divergence-free vector field
$$
u\in L^\infty([0,1];C^\omega(\T^2))
$$
and $f_0\in C^\infty(\T^2)$ such that \eqref{eq:transportintro} admits at least two distinct distributional solutions in $L^\infty([0,1]\times\T^2)$ with initial condition $f_0$.
\end{theorem}

\subsection{Main ideas and comparison with \cite{ElgindiLiss2024}}
	
Following \cite{DEIJ22,ElgindiLiss2024}, we build our velocity field from alternating horizontal and vertical sawtooth (piecewise linear) shear flows whose spatial frequencies blow up geometrically as the switching times accumulate at $t=1$. We conclude the introduction by discussing how the Osgood threshold emerges within such a construction and the main ideas that allow us to go beyond the regularity \eqref{eq:oldregularity} achieved in \cite{ElgindiLiss2024}.

To begin, recall that \emph{any} velocity field that produces anomalous dissipation must generate sufficiently rapid growth of scalar gradients in the absence of diffusion. More precisely, if \eqref{eq:AD} holds at $t=1$, then the solution $f$ of \eqref{eq:transportintro} starting from the same initial data necessarily satisfies
\begin{equation}\label{eq:necessity}
    \lim_{t\uparrow 1}\int_0^t \|\grad f(s)\|_{L^2}^2\,\dee s=\infty.
\end{equation}
This follows from the approximation estimate (see \cite[Section~2.1]{DEIJ22})
\begin{equation}\label{eq:approximation1}
    \|f^\nu(t)-f(t)\|_{L^2}^2
    \le \nu\int_0^t \|\grad f(s)\|_{L^2}^2\,\dee s
\end{equation}
and the conservation of $\|f(t)\|_{L^2}$, which together imply that $\lim_{\nu\to0}\|f^\nu(1)\|_{L^2}=\|f_0\|_{L^2}$ whenever the integral in \eqref{eq:necessity} remains finite.

For alternating shear flows, the key parameter governing gradient growth is the \emph{effective strength} $K_i$ of the $i$th shear, defined as
\[
    K_i=\text{amplitude}\times\text{frequency}\times\text{duration}.
\]
This quantity determines the derivative of the flow map and hence how much the $H^1$ norm of the scalar can be amplified during the $i$th step. A consequence of the arguments in \cite{ElgindiLiss2024} is that the growth required by \eqref{eq:necessity} can be achieved as soon as
\begin{equation}\label{eq:lbd}
    K_i\ge K
\end{equation}
for all $i$ and some sufficiently large $K>0$.

The Osgood threshold appears because it determines precisely when a construction satisfying \eqref{eq:lbd} is possible. Indeed, as we show in Section~\ref{sec:osgood}, fitting a sequence of shears with dyadically increasing frequencies into unit time while maintaining a prescribed modulus of continuity $\omega$ amounts to the requirement that
\begin{equation}\label{eq:sheartimeintro}
    \sum_{i=0}^\infty \frac{K_i\,2^{-i}}{\omega(2^{-i})}<\infty.
\end{equation}
The sum without the factors $K_i$ is comparable to $$\int_0^\pi \frac{\dee r}{\omega(r)},$$ and so the lower bound \eqref{eq:lbd} forces $\omega$ to be non-Osgood. Conversely, for every concave non-Osgood modulus, the constant choice
\begin{equation}\label{eq:constantstrengthintro}
    K_i=K
\end{equation}
satisfies both \eqref{eq:lbd} and the time constraint \eqref{eq:sheartimeintro}. In fact, for any prescribed growth $K_i\to\infty$, there is a non-Osgood modulus for which the sum in \eqref{eq:sheartimeintro} diverges. From the point of view of regularity, at least, a constant shear strength is the natural choice. 

The difficulty is to actually prove anomalous dissipation while keeping the shear strength fixed. The $\nu=0$ gradient growth in \eqref{eq:necessity} is not by itself sufficient, and the restriction on the modulus in \cite{ElgindiLiss2024} arose precisely in passing from this growth to anomalous dissipation. The proof applied a variant of \cite[Proposition~1.3]{DEIJ22}, which provides sufficient conditions for anomalous dissipation in terms of the transport solution alone. In particular, 
\eqref{eq:necessity} does suffice if at the same time $f$ satisfies the ``balanced growth'' condition
\begin{equation}\label{eq:reverseinterpolation}
    \|f(t)\|_{\dot H^{s}}
    \le C\frac{\|f(t)\|_{\dot H^1}^{s}}{\|f(t)\|_{L^2}^{s-1}}
\end{equation}
for some fixed $s>1$ and $C > 0$. Proving \eqref{eq:reverseinterpolation} required the shear strengths $K_i$ to grow at least algebraically and led to the restriction \eqref{eq:oldregularity}.

The role of \eqref{eq:reverseinterpolation} was to supplement the $H^1$ blowup of the $\nu=0$ solution with information about where its energy lies in Fourier space. Specifically, it implies that
\begin{equation}\label{eq:balanced}
    \frac{\|P_{>N(t)}f(t)\|_{L^2}}{\|f(t)\|_{L^2}}\ge c,
    \qquad
    N(t)=\frac{c}{2}\frac{\|\grad f(t)\|_{L^2}}{\|f(t)\|_{L^2}},
\end{equation}
for some $c\in(0,1)$, where $P_{>N}$ denotes the Fourier projection onto frequencies with magnitude greater than $N$. Put simply, a fixed fraction of the scalar Fourier mass remains at frequencies comparable to $N(t)$. It turns out that maintaining $L^2$ mass at this scale is enough to deduce anomalous dissipation from the integrated gradient growth in \eqref{eq:necessity} simply by approximating $f^\nu$ with $f$ using \eqref{eq:approximation1}. While it is quite convenient to reduce anomalous dissipation to norm estimates for the transport solution, balanced growth asks for rather strong control on the Fourier distribution. The mass could easily spread as it moves to higher frequencies, to the point that no fixed fraction keeps pace with $N(t)$. At the same time, such spreading seems entirely consistent with anomalous dissipation, provided the flow still carries a fixed fraction of the Fourier mass of $f$ to arbitrarily high frequencies.

In our construction, we are able to keep the shear strength $K_i=K\gg1$ constant, which is what allows us to reach the full range of non-Osgood moduli in Theorem~\ref{thm:main}. The key to proving anomalous dissipation with this choice is to track the transfer of Fourier mass directly in the $\nu > 0$ equation. Though balanced growth probably fails when $K_i$ remains bounded (see Remark~\ref{rem:balanced}), alternating sawtooth shears still push Fourier mass to increasingly high frequencies as soon as $K$ is sufficiently large. We explain the basic mechanism in Section~\ref{sec:couette} with a discussion of alternating horizontal and vertical Couette flows. These linear shears provide a local model for the sawtooth flows away from their corners and suggest how frequencies should be amplified in our construction.

To justify the Couette picture, we begin in Section~\ref{sec:transferproof} by studying how a sawtooth shear transforms the Fourier distribution of a single mode. We approximate its evolution by gluing together the local Couette solutions and bound the resulting error in Lemma~\ref{lem:approximation}. We call the Fourier mass outside the frequency region predicted by the Couette model the Fourier leakage. Lemma~\ref{lem:profilefourier} bounds this leakage for the approximation. Together, the two estimates give a leakage bound for the original single-mode evolution, stated in Lemma~\ref{lem:onemode}.

To control the leakage for a superposition of modes, the single-mode estimates cannot simply be added because distinct input modes can contribute to the same output frequency. We address this difficulty by introducing a spatial shift in each shear, parametrized by a phase. Averaging over this phase cancels the cross terms between different input modes and yields the phase-averaged leakage bound of Section~\ref{sec:phaseaverage}. In Section~\ref{sec:crossvariation}, we combine this bound with an estimate on how the leakage varies with $\nu$ to choose a single phase at each step that works for every diffusivity, completing the construction of the velocity field. The dependence of our flow on the initial data enters through this choice of phases.

Finally, in Section~\ref{sec:concludeproof}, we control the accumulated Fourier leakage and show that enough Fourier mass is transferred to high frequencies to produce anomalous dissipation, proving Theorem~\ref{thm:main}.

\section{The velocity field and the Osgood condition} \label{sec:osgood}

In this section, we construct the velocity field and relate its parameters to the prescribed modulus of continuity. We follow the framework of alternating sawtooth shear flows in \cite{DEIJ22,ElgindiLiss2024}.

\subsection{Alternating shear flow}
Fix an increasing sequence of times
$$
0=t_0<t_1<t_2<\cdots,\qquad \lim_{n\to\infty}t_n=1,
$$
and denote the duration of the $n$th time step by 
$$\tau_n=t_{n+1}-t_n.$$
Thus $[0,1)$ is partitioned into the intervals $[t_n,t_{n+1})$, with $\sum_{n=0}^{\infty}\tau_n=1$. The durations $\tau_n$ will be specified below.

On each time interval, the velocity field is a sawtooth shear, alternating between the horizontal and vertical directions. Let $S:\R\to\R$ be the sawtooth function, defined as the $2\pi$-periodic extension of the absolute value function from $[-\pi,\pi)$.
The horizontal and vertical shears are denoted by
\begin{align*}
H_{\alpha,N}^\phi(x,y)&=\bigl(\alpha S(Ny+\phi),0\bigr),\\
V_{\alpha,N}^\phi(x,y)&=\bigl(0,\alpha S(Nx+\phi)\bigr).
\end{align*}
The parameters $\alpha>0$, $N\in\N$, and $\phi\in\T$ represent the amplitude, spatial frequency, and phase, respectively. We then set
\begin{equation}\label{eq:u}
 u(t)=u_n,\quad\text{on}\quad [t_n,t_{n+1}),
\end{equation}
where
\begin{equation}\label{eq:un}
u_n=
\begin{cases}
H_{\alpha_n,N_n}^{\phi_n},&n\text{ even},\\[3pt]
V_{\alpha_n,N_n}^{\phi_n},&n\text{ odd}.
\end{cases}
\end{equation}

We will choose the parameters $\tau_n$, $\alpha_n$, $N_n$, and $\phi_n$, independent of $\nu$, so that the velocity field satisfies two requirements:
\begin{itemize}
\item $u$ satisfies a given non-Osgood modulus $\omega$, namely
\begin{equation}\label{eq:uMOC}
 u\in L^\infty([0,1];C^\omega(\T^2)).	
\end{equation}
\item the associated advection-diffusion equation \eqref{eq:ADE} exhibits anomalous dissipation, as stated in Theorem~\ref{thm:main}.
\end{itemize}

\subsection{Non-Osgood modulus}
Given a non-Osgood modulus $\omega$, the requirement \eqref{eq:uMOC} forces a relationship between the amplitudes $\alpha_n$ and frequencies $N_n$. 

\begin{lemma}\label{lem:modulus}
Let $\omega$ be a concave modulus of continuity. For every $n\ge0$ and every phase $\phi_n$, the shear $u_n$ defined in \eqref{eq:un} satisfies
\begin{equation}\label{eq:modulusgoal}
[u_n]_{C^\omega(\T^2)}
\le \frac{\pi\alpha_n}{\omega(\pi/N_n)}.
\end{equation}
\end{lemma}

\begin{proof}
Fix $z,z'\in\T^2$ with $0<d_{\T^2}(z,z')=r\le\pi$. Consider two cases.

When $r\le\pi/N_n$, since $S$ is Lipschitz, we have $|u_n(z)-u_n(z')|\leq \alpha_nN_nr$. Together with the concavity of $\omega$, we have
\[
\frac{|u_n(z)-u_n(z')|}{\omega(r)}
\le\frac{\alpha_nN_nr}{\omega(r)}
\le\frac{\pi\alpha_n}{\omega(\pi/N_n)}.
\]
When $\pi/N_n<r\le\pi$, we use the fact that $S$ takes values in $[0,\pi]$ and the monotonicity of $\omega$ to get
\[
\frac{|u_n(z)-u_n(z')|}{\omega(r)}
\le\frac{\pi\alpha_n}{\omega(r)}
\le\frac{\pi\alpha_n}{\omega(\pi/N_n)}.
\]
\end{proof}

\begin{figure}[htbp]
\centering
\begin{tikzpicture}[x=1.5cm,y=1.15cm,>=Stealth,font=\small]
\draw[->,black!65] (0,0)--(8.4,0) node[right] {$r$};
\draw[->,black!65] (0,0)--(0,3.5);
\node[below left] at (0,0) {$0$};
\draw[densely dashed,black!25] (1,0)--(1,1.1)--(0,1.1);
\draw[densely dashed,black!25] (2,0)--(2,{1.1*sqrt(2)})--(0,{1.1*sqrt(2)});
\draw[densely dashed,black!25] (4,0)--(4,2.2)--(0,2.2);
\node[below] at (1,0) {$\pi/N_{n+2}$};
\node[below] at (2,0) {$\pi/N_{n+1}$};
\node[below] at (4,0) {$\pi/N_n$};
\node[left] at (0,1.1) {$\alpha_{n+2}$};
\node[left] at (0,{1.1*sqrt(2)}) {$\alpha_{n+1}$};
\node[left] at (0,2.2) {$\alpha_n$};
\draw[very thick,orange!45!white]
 (0,0)--(4,2.2)--(8,0);
\draw[very thick,blue!35!white]
 (0,0)--(2,{1.1*sqrt(2)})--(4,0)--(6,{1.1*sqrt(2)})--(8,0);
\draw[very thick,teal!40!white]
 (0,0)--(1,1.1)--(2,0)--(3,1.1)--(4,0)--(5,1.1)--(6,0)--(7,1.1)--(8,0);
\draw[very thick,black!80,domain=0:1,samples=100,variable=\t]
 plot ({8*\t*\t},{1.1*sqrt(8)*\t});
\node[above left] at (8,{1.1*sqrt(8)}) {$c\,\omega(r)$};
\fill[orange!70!black] (4,2.2) circle (1.8pt);
\fill[blue!65!black] (2,{1.1*sqrt(2)}) circle (1.8pt);
\fill[teal!70!black] (1,1.1) circle (1.8pt);
\end{tikzpicture}
\caption{Amplitude and frequency selection. The lighter curves are the normalized sawtooth profiles $\alpha_j S(N_jr)/\pi$, for $j=n,n+1,n+2$, drawn with phase zero. All three start at the origin, and their first peaks touch $c\,\omega(r)$ at $(\pi/N_j,\alpha_j)$. Concavity and monotonicity keep each profile below the modulus curve. As $N_{j+1}=2N_j$, the branch length is halved and the amplitude decreases at each step. The actual shear profiles are $\pi$ times the curves shown.}
\label{fig:modulus-shears}
\end{figure}
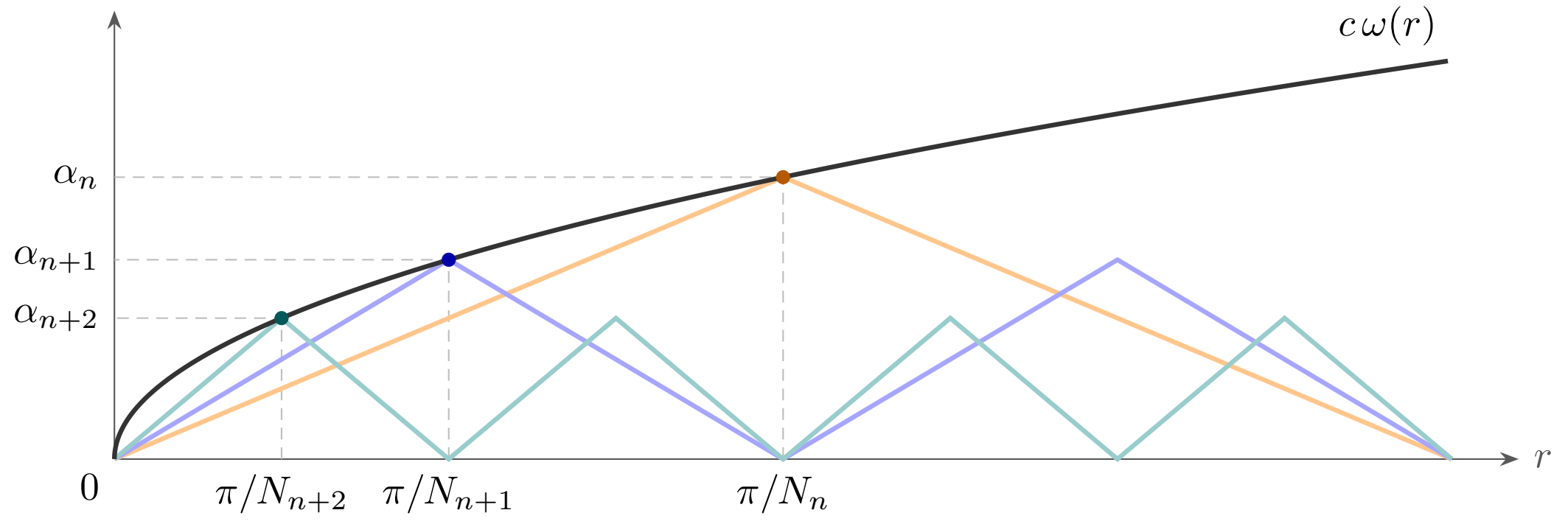

We set $N_n=2^n$. To obtain a uniform $C^\omega$ bound in $n$, we choose the amplitudes
\begin{equation}\label{eq:Nalpharelation}
\alpha_n= c \, \omega(\pi/N_n) = c\,\omega(\pi 2^{-n}),	
\end{equation}
with a constant $c>0$ to be determined. Figure~\ref{fig:modulus-shears} illustrates the relation between the amplitudes and frequencies \eqref{eq:Nalpharelation}.

From \eqref{eq:Nalpharelation}, Lemma~\ref{lem:modulus} implies the uniform bound
\[
\sup_{n\ge0}\,[u_n]_{C^\omega(\T^2)}\le \frac{\pi\alpha_n}{\omega(\pi 2^{-n})} = \pi c.
\]
This leads to the desired modulus of continuity \eqref{eq:uMOC}. Ineed, we have
\[
 \|u\|_{L^\infty([0,1]; C^\omega(\T^2)])} 
 = \sup_{n\ge0} \Big(\|u_n\|_{L^\infty(\T^2)} + [u_n]_{C^\omega(\T^2)}\Big)
 \le \sup_{n\ge0} \Big(\pi\alpha_n + \frac{\pi\alpha_n}{\omega(\pi 2^{-n})}\Big)
 \le \pi c\omega(\pi) + \pi c.
\]

Next, fix $K>32$ and choose the duration $\tau_n$ of each time step so that
\begin{equation}\label{eq:strength}
\alpha_nN_n\tau_n=K
\end{equation}
for every $n\ge0$. This condition fixes the accumulated shear strength at each step and is essential to the mechanism producing anomalous dissipation in our construction: it drives the repeated transfer of Fourier mass to higher frequencies. We explain this mechanism in detail in Section~\ref{sec:couette}.

With $N_n$ and $\alpha_n$ chosen as above, condition \eqref{eq:strength} determines the time steps:
\begin{equation}\label{eq:taundef}
\tau_n=\frac Kc\frac{2^{-n}}{\omega(\pi 2^{-n})}.
\end{equation}
The non-Osgood condition ensures that these time steps have a finite sum. Indeed, monotonicity of $\omega$ gives
\[
\sum_{n=0}^\infty\tau_n
=\frac Kc\sum_{n=0}^\infty\frac{2^{-n}}{\omega(\pi 2^{-n})}
\le\frac{2K}{\pi c}\sum_{n=0}^\infty
\int_{\pi 2^{-(n+1)}}^{\pi 2^{-n}}\frac{\dee r}{\omega(r)}
=\frac{2K}{\pi c}\int_0^\pi\frac{\dee r}{\omega(r)}<\infty,
\]
where the last inequality follows from \eqref{eq:nonOsgood}. We therefore choose
\begin{equation}\label{eq:normalization}
c=K\sum_{n=0}^\infty\frac{2^{-n}}{\omega(\pi 2^{-n})}<\infty,
\end{equation}
so that $\sum_{n=0}^\infty\tau_n=1$, as required. 

\subsection{Summary of the parameter choices}
We choose
\[
K>32,\qquad N_n=2^n,\qquad
\alpha_n=c\,\omega(\pi/N_n),\qquad
\tau_n=\frac{K}{\alpha_nN_n},
\]
with the normalization constant $c$ defined by \eqref{eq:normalization}. These choices ensure that
\begin{itemize}
\item the velocity field has the prescribed non-Osgood regularity \eqref{eq:uMOC};
\item condition \eqref{eq:strength} holds at every step, providing the shear strength needed for the anomalous dissipation mechanism discussed in Section~\ref{sec:couette}.
\end{itemize}
The phases $\phi_n$ will be chosen inductively according to Lemma~\ref{lem:phasechoice}, independently of the diffusivity $\nu$.

\section{Fourier mass transfer for alternating Couette flow}\label{sec:couette}
In this section, we describe the key mechanism that leads to anomalous dissipation for our alternating shear flow. Our analysis tracks the transfer of Fourier mass to successively higher frequencies. 

To start, consider the horizontal Couette flow
\begin{equation}\label{eq:horiCouette}
u(x,y)=\bigl(\alpha(\sigma Ny+\eta),0\bigr),
\qquad \sigma\in\{-1,1\},
\end{equation}
on $\T\times\R$, where $\eta$ is constant. Let $\Phi^{\nu,\mathrm C}_{0,t}$ denote the solution map from time $0$ to time $t$ for \eqref{eq:ADE} with Couette velocity \eqref{eq:horiCouette} and diffusivity $\nu\ge0$. Under pure advection, a plane wave $e_{k,\ell}(x,y)=e^{i(kx+\ell y)}$ evolves exactly as
$$
\big(\Phi^{0,\mathrm C}_{0,t}e_{k,\ell}\big)(x,y)
=e_{k,\ell}\bigl(x-\alpha(\sigma Ny+\eta)t,y\bigr)
=e^{-ik\alpha\eta t}e^{i(kx+(\ell-\sigma\alpha Nk t)y)}.
$$
Thus, after time $\tau$ with \eqref{eq:strength}, namely $\alpha N\tau=K$, the frequency undergoes the shift
$$
(k,\ell)\longmapsto(k,\ell-\sigma Kk).
$$
The constant $\eta$ contributes only a phase factor and does not affect this frequency shift or the amplitude. If $|\ell|\le |k|/2$, then
\begin{equation}\label{eq:FourierCou}
(K-\tfrac12)|k|
\le |\ell-\sigma Kk|
\le (K+\tfrac12)|k|.
\end{equation}
Hence the outgoing $y$-frequency has size comparable to $K|k|$ and dominates the unchanged $x$-frequency.

Diffusion changes the amplitude but not the frequency trajectory. Indeed, substituting
\begin{equation}\label{eq:Couettesolution}
f^\nu(t,x,y)=\big(\Phi^{\nu,\mathrm C}_{0,t}e_{k,\ell}\big)(x,y)=D_\sigma(t)e^{-ik\alpha\eta t}
e^{i(kx+(\ell-\sigma\alpha Nk t)y)}
\end{equation}
into \eqref{eq:ADE} gives
$$
D_\sigma'(t)=-\nu\bigl(k^2+(\ell-\sigma\alpha Nk t)^2\bigr)D_\sigma(t),
\qquad D_\sigma(0)=1.
$$
Consequently,
\begin{equation}\label{eq:Couettedamping}
D_\sigma(t)=\exp\!\left(
-\nu\int_0^t\bigl(k^2+(\ell-\sigma\alpha Nk s)^2\bigr)\,\dee s
\right).
\end{equation}

A vertical Couette flow gives the same mechanism with the roles of $k$ and $\ell$ reversed. More precisely, the vertical Couette flow $u(x,y)=(0,\alpha(\sigma Nx+\eta))$ maps the frequency $(k,\ell)$ at time $\tau$, with $\alpha N\tau=K$, to
$$
(k,\ell)\longmapsto(k-\sigma K\ell,\ell).
$$
If $|k|\le|\ell|/2$, then
$$
(K-\tfrac12)|\ell|
\le |k-\sigma K\ell|
\le (K+\tfrac12)|\ell|,
$$
so the same amplification occurs, with the $x$-frequency now dominant. As in the horizontal case, diffusion changes only the amplitude. Alternating the two orientations therefore moves frequencies between horizontal and vertical cones, amplifying their dominant component by a factor of roughly $K$ at each step.

We now alternate horizontal and vertical Couette flows, using the parameters $\alpha_n,N_n,\tau_n$ with $\alpha_nN_n\tau_n=K$. On the $n$th time interval $[t_n,t_{n+1}]$, consider the affine velocity
$$
u_n^{\mathrm C}(x,y)=
\begin{cases}
\bigl(\alpha_n(\sigma_nN_ny+\eta_n),0\bigr),&n\text{ even},\\
\bigl(0,\alpha_n(\sigma_nN_nx+\eta_n)\bigr),&n\text{ odd},
\end{cases}
\qquad \sigma_n\in\{-1,1\},
$$
where $\eta_n$ is constant. These affine flows are considered on $\R^2$. For this Couette model, we take $K>32$ to be an integer so that integer frequencies remain integer at the stage endpoints.

Starting from the initial scalar $f_0(x,y)=\cos(mx)$, define
\begin{equation}\label{eq:Couettecones}
C_n=
\begin{cases}
\bigl\{(k,\ell)\in\Z^2:
m(K/2)^n\le|k|\le m(2K)^n,\
|\ell|\le|k|/2\bigr\},
&n\text{ even},\\[4pt]
\bigl\{(k,\ell)\in\Z^2:
m(K/2)^n\le|\ell|\le m(2K)^n,\
|k|\le|\ell|/2\bigr\},
&n\text{ odd}.
\end{cases}
\end{equation}
If $(k,\ell)\in C_n$ and $n$ is even, the outgoing frequency is $(k,\ell-\sigma_nKk)$. By \eqref{eq:FourierCou},
$$
m(K/2)^{n+1} \le (K-\tfrac12)m(K/2)^n \le |\ell-\sigma_nKk| \le (K+\tfrac12)m(2K)^n \le m(2K)^{n+1},
$$
and
$$
|k|\le \frac{|\ell-\sigma_nKk|}{K-\tfrac12}
\le \tfrac12|\ell-\sigma_nKk|.
$$
Thus $(k,\ell-\sigma_nKk)\in C_{n+1}$. For odd $n$, interchanging the coordinates gives $(k-\sigma_nK\ell,\ell)\in C_{n+1}$. The frequency transfer from $C_n$ to $C_{n+1}$ therefore holds at every stage, independently of the sign $\sigma_n$ and the constant $\eta_n$.

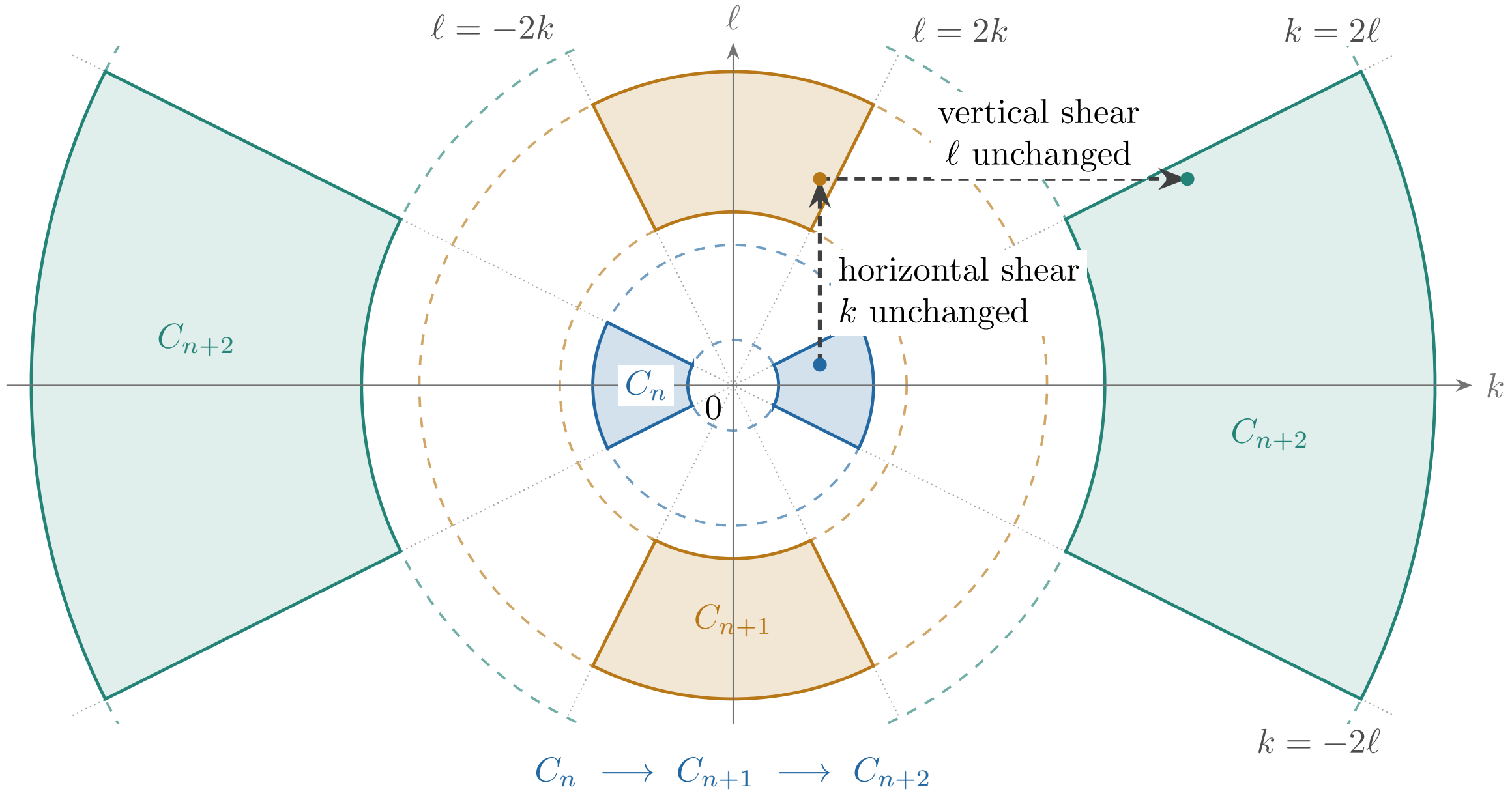
\begin{figure}[htbp]
\centering
\begin{tikzpicture}[x=0.33in,y=0.33in,
  >=Stealth, font=\small,
  cone/.style={line width=0.9pt},
  transfer/.style={->,very thick,densely dashed,black!75}]
\definecolor{coneone}{RGB}{35,104,162}
\definecolor{conetwo}{RGB}{184,120,23}
\definecolor{conethree}{RGB}{35,131,118}
\begin{scope}
\clip (-8.8,-4.1) rectangle (8.8,4.1);
\foreach \radius in {0.55,1.7}
  \draw[overlay,draw=coneone!65,dashed,line width=0.7pt] (0,0) circle[radius=\radius];
\foreach \radius in {2.1,3.8}
  \draw[overlay,draw=conetwo!65,dashed,line width=0.7pt] (0,0) circle[radius=\radius];
\foreach \radius in {4.5,8.5}
  \draw[overlay,draw=conethree!65,dashed,line width=0.7pt] (0,0) circle[radius=\radius];
\end{scope}
\draw[black!35,densely dotted] (-2,-4)--(2,4)
  node[above right,text=black!70] {$\ell=2k$};
\draw[black!35,densely dotted] (2,-4)--(-2,4)
  node[above left,text=black!70] {$\ell=-2k$};
\draw[black!35,densely dotted] (-8,-4)--(8,4)
  node[above left,text=black!70] {$k=2\ell$};
\draw[black!35,densely dotted] (-8,4)--(8,-4)
  node[below left,text=black!70] {$k=-2\ell$};
\pgfmathsetmacro{\coneangle}{atan(0.5)}
\foreach \rotation in {0,180} {
  \begin{scope}[rotate=\rotation]
    \draw[cone,draw=coneone,fill=coneone!20]
      (-\coneangle:0.55)--(-\coneangle:1.7)
      arc[start angle=-\coneangle,end angle=\coneangle,radius=1.7]
      --(\coneangle:0.55)
      arc[start angle=\coneangle,end angle=-\coneangle,radius=0.55]
      --cycle;
    \draw[cone,draw=conethree,fill=conethree!14]
      (-\coneangle:4.5)--(-\coneangle:8.5)
      arc[start angle=-\coneangle,end angle=\coneangle,radius=8.5]
      --(\coneangle:4.5)
      arc[start angle=\coneangle,end angle=-\coneangle,radius=4.5]
      --cycle;
  \end{scope}
}
\foreach \rotation in {90,270} {
  \begin{scope}[rotate=\rotation]
    \draw[cone,draw=conetwo,fill=conetwo!18]
      (-\coneangle:2.1)--(-\coneangle:3.8)
      arc[start angle=-\coneangle,end angle=\coneangle,radius=3.8]
      --(\coneangle:2.1)
      arc[start angle=\coneangle,end angle=-\coneangle,radius=2.1]
      --cycle;
  \end{scope}
}
\draw[->,black!55] (-8.8,0)--(8.95,0) node[right] {$k$};
\draw[->,black!55] (0,-4.1)--(0,4.15) node[above] {$\ell$};
\node[below left,inner sep=3pt] at (0,0) {$0$};
\node[text=coneone,fill=white,inner sep=1.7pt] at (-1.05,0) {$C_n$};
\node[text=conetwo] at (0,-2.85) {$C_{n+1}$};
\node[text=conethree] at (-6.5,0.55) {$C_{n+2}$};
\node[text=conethree] at (6.5,-0.6) {$C_{n+2}$};
\draw[transfer] (1.05,0.25)--(1.05,2.5);
\node[anchor=west,fill=white,inner sep=2pt,align=left]
  at (1.18,1.12) {horizontal shear\\$k$ unchanged};
\draw[transfer] (1.05,2.5)--(5.5,2.5);
\node[above,fill=white,inner sep=2pt,align=center] at (3.7,2.5)
  {vertical shear\\$\ell$ unchanged};
\fill[coneone] (1.05,0.25) circle[radius=2pt];
\fill[conetwo] (1.05,2.5) circle[radius=2pt];
\fill[conethree] (5.5,2.5) circle[radius=2pt];
\node[below,text=coneone] at (0,-4.35)
  {$C_n\ \longrightarrow\ C_{n+1}\ \longrightarrow\ C_{n+2}$};
\end{tikzpicture}
\caption{Schematic illustration of $C_n$, $C_{n+1}$, and $C_{n+2}$ for even $n$, drawn as sectors of successive concentric rings to emphasize transfer to higher frequencies. The vertical and horizontal arrows emphasize preservation of $k$ and $\ell$, respectively, at successive shear steps. The dotted lines $\ell=\pm2k$ and $k=\pm2\ell$ give the exact angular boundaries of the vertical and horizontal cones. Each ring has its own inner and outer circle. The gaps and radial scales are schematic and do not impose a separation condition on the cones; the precise coordinate cutoffs are those in \eqref{eq:Couettecones}.}
\label{fig:alternating-cones}
\end{figure}

Figure~\ref{fig:alternating-cones} illustrates the alternating orientations and the increasing frequency scales of these cones.

Let $P_n$ denote the Fourier projection onto $C_n$, so that $P_0f_0=f_0$, and let $\Phi^{\nu,\mathrm C}_n=\Phi^{\nu,\mathrm C}_{t_n,t_{n+1}}$ denote the solution map from time $t_n$ to time $t_{n+1}$ for the Couette velocity $u_n^{\mathrm C}$. For $g\in L^2(\T^2)$ with Fourier support in $C_n$, we apply this map to the periodic extension of $g$ to $\R^2$. Since $K$ is an integer, the solution at time $t_{n+1}$ is again periodic. The frequency transfer established above gives
\begin{equation}\label{eq:Couetteproj}
\bigl\|(1-P_{n+1})\Phi^{\nu,\mathrm C}_n g\bigr\|_{L^2}^2=0,
\qquad g=P_n g,\quad n\ge0,\quad \nu\ge0.
\end{equation}
This means that all Fourier mass of the evolved scalar lies in $C_{n+1}$ whenever the incoming scalar has Fourier support in $C_n$. This holds for every choice of $\sigma_n$ and $\eta_n$ and every $\nu\ge0$, because diffusion changes only the amplitudes of the evolving modes.

Iterating \eqref{eq:Couetteproj} confines the Fourier support at time $t_n$ to $C_n$. Since every frequency in $C_n$ has magnitude at least $m(K/2)^n\to\infty$, this support leaves every bounded frequency region as $t_n\to1$. In the absence of diffusion, the total $L^2$ mass at these times is conserved, so the mass escapes to infinity in frequency space. This transfer to increasingly fine scales is the mechanism underlying anomalous dissipation: even weak diffusivity can dissipate mass transported to sufficiently high frequencies.

Our alternating shear flow \eqref{eq:un} agrees with a Couette flow on each linear branch, suggesting a similar mechanism of Fourier mass transfer. However, the kinks in its profile spread frequencies even when $\nu=0$, and diffusion couples the evolution of the scalar in neighboring branches when $\nu>0$. Consequently, the exact confinement to $C_{n+1}$ in \eqref{eq:Couetteproj} need not hold. We refer to the amount of Fourier mass outside $C_{n+1}$ after this step as \emph{Fourier leakage}.

The following lemma quantifies the Fourier leakage of each iteration. Denote $\Psi_{n,\phi}^\nu$ the solution map from time $t_n$ to time $t_{n+1}$ for \eqref{eq:ADE} with $u= u_n$ and $\phi_n=\phi$.

\begin{lemma}\label{lem:transfer}
There exists $C > 0$ such that for every $n\ge0$, $\nu\ge0$, $K > 32$, and $g \in L^2(\T^2)$ with Fourier support in $C_n$, we have
\begin{equation}\label{eq:transfer}
\frac1{2\pi}\int_0^{2\pi}
\bigl\|(1-P_{n+1})\Psi_{n,\phi}^\nu g\bigr\|_{L^2}^2 \, \dee\phi
\le \frac{C}{m}\left(\frac4K\right)^n\|g\|_{L^2}^2.
\end{equation}
\end{lemma}

We prove Lemma~\ref{lem:transfer} in Section~\ref{sec:transferproof}. The proof is based on first establishing the estimate without the phase averaging when $g = e_{k,\ell}$ for $(k,\ell) \in C_n$. Phase averaging then allows us to sum the single-mode estimate to obtain the result for general scalars by eliminating the cross terms that arise from interactions between distinct incoming frequencies. 

To track the transfer of Fourier mass with Lemma~\ref{lem:transfer}, we introduce an auxiliary sequence of scalars that is projected onto the relevant cone at each step. For phases $\{\phi_n\}_{n \ge 0}$ to be fixed below, set $g_0^\nu=f_0$ and define recursively
\begin{equation}\label{eq:iteration}
g_{n+1}^\nu=P_{n+1}\Psi_{n,\phi_n}^\nu g_n^\nu.
\end{equation}
Since $P_0f_0=f_0$, each $g_n^\nu$ has Fourier support in $C_n$ 
and it follows from Lemma~\ref{lem:transfer} that for each fixed $\nu > 0$ there exists a choice of $\phi_n$ such that
$$
\bigl\|(1-P_{n+1})\Psi_{n,\phi_n}^\nu g_n^\nu\bigr\|_{L^2}^2
\le
\frac C m\left(\frac4K\right)^n\|g_n^\nu\|_{L^2}^2.
$$
However, this phase may depend on $\nu$ through both the solution operator and the incoming scalar $g_n^\nu$. To obtain a $\nu$-independent velocity field, we must show that one phase works for any small diffusivity. We achieve this by combining Lemma~\ref{lem:transfer} with quantitative continuity-in-$\nu$ estimates for the family $\{g_n^\nu\}_{\nu\ge0}$. Specifically, in Section~\ref{sec:crossvariation} we prove the following.

\begin{lemma}\label{lem:phasechoice}
Fix $K>32$. There exist $m\in\N$ and $\nu$-independent phases $\{\phi_n\}_{n\ge0}$ such that, for every $n\ge0$,
\begin{equation}\label{eq:phasechoice}
\sup_{\nu\in[0,1]}
\bigl\|(1-P_{n+1})\Psi_{n,\phi_n}^\nu g_n^\nu\bigr\|_{L^2}
\le 2^{-n-2}\|f_0\|_{L^2}.
\end{equation}
\end{lemma}

We apply Lemma~\ref{lem:phasechoice} to choose phases independent of the diffusivity and bound the cumulative error caused by Fourier leakage. The summable bounds keep this error smaller than the initial $L^2$ norm, so that, in the absence of diffusion, a positive amount of Fourier mass remains in the successive cones $C_n$. As in the Couette case, their minimum frequency $m(K/2)^n$ tends to infinity as $t_n\to1$, and this mass is pushed to arbitrarily high frequencies despite the leakage. For $\nu > 0$, the projected scalar dissipates completely, while the same error bound yields a positive lower bound on the dissipation of the actual solution, uniformly as $\nu\to0$. This proves Theorem~\ref{thm:main}. The detailed argument is given in Section~\ref{sec:concludeproof}.

\begin{remark} \label{rem:balanced}
Tracking Fourier mass through the successive cones $C_n$ is sufficient for us to prove anomalous dissipation, but it does not provide enough information to deduce that the norm growth is balanced in the sense of \eqref{eq:balanced}. Controlling the leakage in our argument requires us to widen the cones at each step. The frequencies in $C_n$ can then differ in size by a factor of $4^n$, and with this much room for the Fourier mass to spread, it is easy to see that \eqref{eq:balanced} may in principle be false. Since even a single sawtooth shear spreads an incoming Fourier mode over a range of frequencies, we expect some of this widening to reflect the actual dynamics, although strictly speaking we do not prove that \eqref{eq:balanced} fails for our solution. 
\end{remark}

\section{Proof of the Fourier transfer lemma}\label{sec:transferproof}

In this section, we prove the key Fourier leakage bound in Lemma~\ref{lem:transfer}. We first consider a single input Fourier mode and establish the following estimate.

\begin{lemma}\label{lem:onemode}
    There exists a constant $C > 0$ independent of $\alpha_n$, $N_n$, and $\tau_n$ such that for every even $n\ge0$, any $(k,\ell) \in C_n$, phase $\phi$, and $\nu \ge 0$, we have 
\begin{equation}\label{eq:onemode}
\bigl\|(1-P_{n+1})\Psi_{n,\phi}^\nu e_{k,\ell}\bigr\|_{L^2}^2
\le C\frac{N_n}{|k|}.
\end{equation}
\end{lemma}

We then expand a general incoming scalar in Fourier modes and use phase averaging to combine the single-mode estimates. The averaging eliminates cross terms between distinct incoming modes, yielding Lemma~\ref{lem:transfer}.

Throughout this section, we fix even $n\ge0$, so that $u_n$ is a horizontal shear, and write $\alpha=\alpha_n$, $N=N_n$, and $\tau=\tau_n$. The case of odd $n$ follows by the same arguments with the coordinates interchanged. Recall that $C_n$ is defined by \eqref{eq:Couettecones}, $P_n$ is the Fourier projection onto $C_n$, and $\Psi_{n,\phi}^\nu$ is the time-$\tau_n$ solution map for $u_n$ with phase $\phi_n=\phi$.

\subsection{Approximation of $\Psi_{n,\phi}^\nu e_{k,\ell}$}\label{sec:singlemode}

The proof of Lemma~\ref{lem:onemode} proceeds as follows. We first construct an approximation of $\Psi_{n,\phi}^\nu e_{k,\ell}$ by solving \eqref{eq:ADE} separately on each region where the shear is linear and then piecing these solutions together. This ignores the coupling between neighboring regions caused by diffusion. On each region, the approximate solution experiences the frequency shift and damping described by \eqref{eq:Couettesolution}. We show that the approximation error is consistent with \eqref{eq:onemode} and then estimate the Fourier mass of the approximation outside $C_{n+1}$ by a direct computation.

Fix $(k,\ell)\in C_n$ and a phase $\phi$. The solution of \eqref{eq:ADE} with $u = u_n = H_{\alpha,N}^\phi$ and initial datum $e_{k,\ell}$ can be written as $e^{ikx}h(t,y)$, where
\begin{equation}\label{eq:reduced}
\begin{cases}
\partial_t h+ik\alpha S(Ny+\phi)h=\nu \partial_{yy} h-\nu k^2h,\\
h(0,y)=e^{i\ell y}.
\end{cases}
\end{equation}

To define an approximate solution of \eqref{eq:reduced}, partition $\T$ into intervals $\{I_j\}_{j=1}^{2N}$ on which $S(Ny+\phi)$ is linear. On each $I_j$, we can write
$$
S(Ny+\phi)=\sigma_jNy+\eta_j,
\qquad \sigma_j\in\{-1,1\},
$$
for some constant $\eta_j$. Now extend this linear shear to $\R$ and let $H_j:[0,\tau] \times \R \to \mathbb{C}$ solve
\begin{equation}\label{eq:linearbranch}
\begin{cases}
\partial_t H_j+ik\alpha(\sigma_jNy+\eta_j)H_j
=\nu\partial_{yy}H_j-\nu k^2H_j,\\
H_j(0,y)=e^{i\ell y}.
\end{cases}
\end{equation}
Applying the Couette formula \eqref{eq:Couettesolution} with $\sigma=\sigma_j$ and $\eta=\eta_j$, we obtain
\begin{equation}\label{eq:couette}
H_j(t,y)= e^{-ik\alpha\eta_j t} D_{\sigma_j}(t) e^{i(\ell-\sigma_j\alpha Nk t)y},
\end{equation}
where $D_\sigma$ is defined by \eqref{eq:Couettedamping}. We piece these solutions together to define our approximation. That is, we set
\begin{equation}\label{eq:branch}
H(t,y)=H_j(t,y),
\qquad y\in I_j.
\end{equation}
By construction, $H$ solves \eqref{eq:reduced} on each $I_j$ and agrees with the initial data. Hence, without diffusion, $H_j=h$ on $I_j$. For $\nu>0$, however, diffusion couples neighboring regions, so $H_j$ need not coincide with $h$ on $I_j$. The following lemma controls the $L^2$ approximation error uniformly in the diffusivity and phase.

\begin{lemma}\label{lem:approximation}
There exists $C>0$ such that for every $\nu\ge0$ and phase $\phi$, we have
\begin{equation}\label{eq:approximation}
\|h(\tau)-H(\tau)\|_{L^2}^2
\le C\frac N{|k|}.
\end{equation}
\end{lemma}

\begin{proof}
The standard $L^2$ energy estimate for $h$ and the explicit formula \eqref{eq:couette} for $H$ give $\|h(\tau)\|_{L^2}\le1$ and $\|H(\tau)\|_{L^2}\le1$. Hence, we have the bound
\[
\|h(\tau)-H(\tau)\|_{L^2}\le\|h(\tau)\|_{L^2}+\|H(\tau)\|_{L^2}\le2,
\]
so \eqref{eq:approximation} is trivial when $|k|\leqc N$.

We therefore focus on the case $|k|\gg N$. Set $\delta=|k|^{-1}\ll N^{-1}$. Let $\chi:\T\to[0,1]$ be a smooth cutoff that is zero within distance $\delta$ of every kink in $S(Ny+\phi)$ and one at distance at least $2\delta$ from every kink. We choose $\chi$ with $|\chi'|\leqc\delta^{-1}$ and $|\chi''|\leqc\delta^{-2}$. Since there are $2N$ kinks,
\begin{equation}\label{eq:measurebound}
|\{y\in\T:\chi(y)\neq1\}|\leqc N\delta.
\end{equation}
This gives the bounds
\begin{equation}\label{eq:cutoffbounds}
\|1-\chi\|_{L^2}\leqc\sqrt{N\delta},\qquad
\|\chi'\|_{L^2}\leqc\delta^{-1}\sqrt{N\delta},\qquad
\|\chi''\|_{L^2}\leqc\delta^{-2}\sqrt{N\delta}.
\end{equation}

We decompose $H=\chi H+(1-\chi)H$ into an interior part and a boundary layer. Then
\[
h-H=(h-\chi H)-(1-\chi)H.
\]

To control the boundary layer, note that the explicit formula \eqref{eq:couette} gives
\begin{equation}\label{eq:branchbound}
\|H(t)\|_{L^\infty}\le e^{-\nu k^2t},
\qquad 0\le t\le\tau.
\end{equation}
Together with \eqref{eq:cutoffbounds}, this yields
\begin{equation}\label{eq:boundarylayer}
\|(1-\chi)H(\tau)\|_{L^2}
\le\|1-\chi\|_{L^2}\|H(\tau)\|_{L^\infty}
\leqc \sqrt{N\delta}.
\end{equation}

For the interior part, $H$ solves the same equation \eqref{eq:reduced} as $h$ on the support of $\chi$. Thus,
\begin{equation}\label{eq:chiH}
\partial_t(\chi H)+ik\alpha S(Ny+\phi)(\chi H)
=\nu\partial_{yy}(\chi H)-\nu k^2(\chi H)
-2\nu\partial_y(\chi'H)+\nu\chi''H
\end{equation}
holds for every $y\in\T$. Let
\[
w=h-\chi H.
\]
Then, from \eqref{eq:reduced} and \eqref{eq:chiH} we see that $w$ solves
\begin{equation}\label{eq:cutoffequation}
w_t+ik\alpha S(Ny+\phi)w
=\nu\partial_{yy}w-\nu k^2w
+2\nu\partial_y(\chi'H)-\nu\chi''H,
\end{equation}
with
\[
w(0,y)=(1-\chi(y))e^{i\ell y}.
\]
A straightforward $L^2$ energy estimate for \eqref{eq:cutoffequation}, followed by Young's inequality, gives
\begin{equation}\label{eq:L2Young}
\frac12\frac{\dee}{\dee t}\|w\|_{L^2}^2
+\frac\nu2\|\partial_yw\|_{L^2}^2
+\frac{\nu k^2}{2}\|w\|_{L^2}^2
\le2\nu\|\chi'H\|_{L^2}^2
+\frac{\nu}{2k^2}\|\chi''H\|_{L^2}^2.
\end{equation}
Putting \eqref{eq:cutoffbounds} and \eqref{eq:branchbound} into \eqref{eq:L2Young}, and using $\delta=|k|^{-1}$, yields
\begin{equation}\label{eq:L2Young2}
\frac12\frac{\dee}{\dee t}\|w\|_{L^2}^2
\leqc\nu e^{-2\nu k^2t}
\left(\frac N\delta+\frac N{k^2\delta^3}\right)
\leqc\nu k^2e^{-2\nu k^2t}\frac N{|k|}.
\end{equation}
Integrating over $t\in[0,\tau]$ gives
\begin{equation}\label{eq:wtau}
\|w(\tau)\|_{L^2}^2
\leqc\|w(0)\|_{L^2}^2+\frac N{|k|}
=\|1-\chi\|_{L^2}^2+\frac N{|k|}
\leqc\frac N{|k|},
\end{equation}
where in the final inequality we used \eqref{eq:cutoffbounds}. 

Combining \eqref{eq:boundarylayer} and \eqref{eq:wtau}, and using the triangle inequality with $\delta=|k|^{-1}$, we conclude with the desired bound
\begin{equation}\label{eq:wtriangle}
\|h(\tau)-H(\tau)\|_{L^2}^2
\le2\|w(\tau)\|_{L^2}^2+2\|(1-\chi)H(\tau)\|_{L^2}^2
\leqc\frac N{|k|}.
\end{equation}
\end{proof}

\subsection{Fourier leakage of the approximation}

In view of \eqref{eq:approximation}, to prove Lemma~\ref{lem:onemode} it is sufficient to establish the estimate for our approximate solution. Since $\alpha N\tau=K$, formula \eqref{eq:couette} gives
\begin{equation}\label{eq:approxwave}
H(\tau,y)=\sum_{j=1}^{2N}
D_{\sigma_j}(\tau)e^{-ik\alpha\eta_j\tau}
e^{i(\ell-\sigma_jKk)y}\mathbf{1}_{I_j}(y).
\end{equation}
The following lemma controls the Fourier mass of $e^{ikx}H(\tau,y)$ outside of $C_{n+1}$.

\begin{lemma}\label{lem:profilefourier}
There exists $C>0$, independent of the diffusivity and phase, such that
\begin{equation}\label{eq:approxtransfer}
\bigl\|(1-P_{n+1})e^{ikx}H(\tau,y)\bigr\|_{L^2}^2
\le C\frac N{|k|}.
\end{equation}
\end{lemma}

\begin{proof}
From \eqref{eq:approxwave} and $\alpha N\tau=K$, we have
\[
e^{-i\ell y}H(\tau,y)
=\sum_{j=1}^{2N}D_{\sigma_j}(\tau)
e^{-i\frac{Kk}{N}(\sigma_jNy+\eta_j)}\mathbf{1}_{I_j}(y).
\]
The piecewise slope $\sigma_j$ and affine profile $\sigma_jNy+\eta_j$ repeat with period $2\pi/N$, so the right-hand side is $2\pi/N$-periodic. Consequently, $e^{ikx}H(\tau,y)$ contains only Fourier modes of the form $(k,\ell+pN)$, for $p\in\Z$.

Writing $\widehat H(\tau,q)$ for the Fourier coefficient of $H(\tau,\cdot)$ in $y$, we obtain from \eqref{eq:approxwave}
\begin{equation}\label{eq:profileformula}
\widehat H(\tau,\ell+pN)=\frac1{2\pi}
\sum_{j=1}^{2N}D_{\sigma_j}(\tau)e^{-ik\alpha\eta_j\tau}
\int_{I_j}e^{-i(Np+\sigma_jKk)y}\,\dee y.
\end{equation}
Since $|I_j|=\pi/N$, each integral above satisfies
\[
\left|\int_{I_j}e^{-i(Np+\sigma_jKk)y}\,\dee y\right|
\le\min\left\{\frac\pi N,\frac2{|Np+\sigma_jKk|}\right\}
\leqc\frac1{N\bigl(1+|p+\sigma_jKk/N|\bigr)}.
\]
Putting this bound into \eqref{eq:profileformula} and using $0\le D_{\sigma_j}(\tau)\le1$ yields
\begin{equation}\label{eq:profilecoefficient}
|\widehat H(\tau,\ell+pN)|
\leqc
\frac1{1+|p-Kk/N|}
+\frac1{1+|p+Kk/N|}.
\end{equation}

We now use \eqref{eq:profilecoefficient} to quantify how much mass lies away from the two expected frequency shifts. First observe that for every $b\in\R$ and $R\ge0$,
$$
\sum_{\substack{p\in\Z\\|p-b|>R}}
\frac1{(1+|p-b|)^2}
\leqc\int_{|x|\ge R}\frac1{1+x^2}\,\dee x
\leqc\frac1{1+R}.
$$
Squaring \eqref{eq:profilecoefficient} and applying this estimate with $b=\pm Kk/N$ and $R=K|k|/(4N)$ gives
\begin{equation}\label{eq:tail}
\sum_{\substack{p\in\Z\\
\min\{|Np-Kk|,|Np+Kk|\}>K|k|/4}}
|\widehat H(\tau,\ell+pN)|^2
\leqc\frac {1}{1+K|k|/N}
\leqc \frac N{|k|}.
\end{equation}
By \eqref{eq:tail} and Plancherel, to complete the proof, it suffices to show that 
\begin{equation}\label{eq:containment}
\min\{|Np-Kk|,|Np+Kk|\}\le\frac{K|k|}{4} \implies (k,\ell+pN) \in C_{n+1}.
\end{equation}
Suppose now that the left-hand side of \eqref{eq:containment} holds. Since $(k,\ell) \in C_n$ implies $|\ell|\le|k|/2$, the triangle inequality gives
$$
\left(\frac{3K}{4}-\frac12\right)|k|
\le|\ell+pN|
\le\left(\frac{5K}{4}+\frac12\right)|k|.
$$
Since  $m(K/2)^n \le |k| \le m(2K)^n$ and $K > 32$, it follows that
\begin{equation}
m(K/2)^{n+1}\le|\ell+pN|\le m(2K)^{n+1}
\end{equation}
and
$$
|k|\le\frac12|\ell+pN|.
$$
The previous two inequalities establish \eqref{eq:containment} and complete the proof.
\end{proof}

We are now ready to conclude the proof of Lemma~\ref{lem:onemode} by combining Lemmas~\ref{lem:approximation} and~\ref{lem:profilefourier}.

\begin{proof}[Proof of Lemma~\ref{lem:onemode}]
Since $\Psi_{n,\phi}^\nu e_{k,\ell}=e^{ikx}h(\tau,y)$, the triangle inequality followed by Lemmas~\ref{lem:approximation} and~\ref{lem:profilefourier} gives
\begin{align*}
\bigl\|(1-P_{n+1})\Psi_{n,\phi}^\nu e_{k,\ell}\bigr\|_{L^2}^2 \le 2\bigl\|(1-P_{n+1})e^{ikx}H(\tau,y)\bigr\|_{L^2}^2
+2\bigl\|e^{ikx}(h(\tau)-H(\tau))\bigr\|_{L^2}^2 \leqc \frac N{|k|}.
\end{align*}
\end{proof}

\subsection{Phase averaging}\label{sec:phaseaverage}

We now upgrade the single-mode estimate of Lemma~\ref{lem:onemode} to a general incoming scalar. Fix $\nu\ge0$ and $g=P_ng$. Expanding $g$ in Fourier modes gives
\[
g=\sum_{(k,\ell)\in C_n}\widehat g(k,\ell)e_{k,\ell},
\qquad
\Psi_{n,\phi}^\nu g
=\sum_{(k,\ell)\in C_n}\widehat g(k,\ell)\Psi_{n,\phi}^\nu e_{k,\ell}.
\]
Since horizontal shears preserve $x$-frequencies, only input modes with the same $k$ contribute to an output frequency $(k,q)$. Consequently, Plancherel gives
\begin{equation}\label{eq:leakageexpansion}
\bigl\|(1-P_{n+1})\Psi_{n,\phi}^\nu g\bigr\|_{L^2}^2
=\sum_{(k,q)\notin C_{n+1}}
\left|\sum_{\ell:(k,\ell)\in C_n}
\widehat g(k,\ell)\,
\widehat{\Psi_{n,\phi}^\nu e_{k,\ell}}(k,q)\right|^2
=\mathcal D_\phi+\mathcal C_\phi.
\end{equation}
Here, the diagonal terms are
\begin{align*}
\mathcal D_\phi
&=\sum_{(k,q)\notin C_{n+1}}
\sum_{\ell:(k,\ell)\in C_n}
|\widehat g(k,\ell)|^2
\left|\widehat{\Psi_{n,\phi}^\nu e_{k,\ell}}(k,q)\right|^2\\
&=\sum_{(k,\ell)\in C_n}|\widehat g(k,\ell)|^2
\bigl\|(1-P_{n+1})\Psi_{n,\phi}^\nu e_{k,\ell}\bigr\|_{L^2}^2,
\end{align*}
while the cross terms are
\begin{equation}\label{eq:crossterms}
\mathcal C_\phi
=2\operatorname{Re}\sum_{(k,q)\notin C_{n+1}}
\sum_{\substack{\ell<\ell'\\(k,\ell),(k,\ell')\in C_n}}
\widehat g(k,\ell)\overline{\widehat g(k,\ell')}
\widehat{\Psi_{n,\phi}^\nu e_{k,\ell}}(k,q)
\overline{\widehat{\Psi_{n,\phi}^\nu e_{k,\ell'}}(k,q)}.
\end{equation}
Applying Lemma~\ref{lem:onemode} directly to the diagonal terms gives
\begin{equation}\label{eq:diagonalbound}
\mathcal D_\phi
\leqc N\sum_{(k,\ell)\in C_n}\frac{|\widehat g(k,\ell)|^2}{|k|}
\leqc\frac{N}{m(K/2)^n}\|g\|_{L^2}^2
=\frac1m\left(\frac4K\right)^n\|g\|_{L^2}^2,
\end{equation}
where we used $N = 2^n$ and the fact that $|k|\ge m(K/2)^n$ on $C_n$.

The cross terms $\mathcal C_\phi$ arise because distinct incoming modes can contribute to the same outgoing frequency. The single-mode estimate alone does not control these terms sufficiently: taking absolute values and applying Cauchy--Schwarz introduces a factor depending on the number of incoming modes, which can destroy the smallness in \eqref{eq:diagonalbound}.

To address this issue, we introduce a phase averaging argument that eliminates the cross terms, as stated in the following lemma.

\begin{lemma}\label{lem:crossaverage}
For $\mathcal C_\phi$ defined in \eqref{eq:crossterms}, we have
\[
\frac1{2\pi}\int_0^{2\pi}\mathcal C_\phi\,\dee\phi=0.
\]
\end{lemma}

\begin{proof}
Using that derivatives commute with translations, one can check that for general $f\in L^2$ there holds
\begin{equation}\label{eq:generaltranslation}
[\Psi_{n,\phi}^\nu f](x,y)
=[\Psi_{n,0}^\nu F](x,y+\phi/N),
\qquad F(x,y)=f(x,y-\phi/N).
\end{equation}
Applying this with $f=e_{k,\ell}$ gives
\begin{equation}\label{eq:modetranslation}
\Psi_{n,\phi}^\nu e_{k,\ell}(x,y)
=e^{-i\ell\phi/N}
\bigl(\Psi_{n,0}^\nu e_{k,\ell}\bigr)(x,y+\phi/N).
\end{equation}
Since the phases $0$ and $2\pi$ define the same shear, setting $\phi=2\pi$ and multiplying the result by $e^{-i\ell y}$ yields
\[
e^{-i\ell y}\bigl(\Psi_{n,0}^\nu e_{k,\ell}\bigr)(x,y)
=e^{-i\ell(y+2\pi/N)}\bigl(\Psi_{n,0}^\nu e_{k,\ell}\bigr)(x,y+2\pi/N).
\]
Thus, $y\mapsto e^{-i\ell y}\Psi_{n,0}^\nu e_{k,\ell}(x,y)$ is $2\pi/N$-periodic. Together with preservation of $x$-frequencies, this shows that $\Psi_{n,0}^\nu e_{k,\ell}$ contains only modes of the form $(k,\ell+pN)$, for $p\in\Z$. By \eqref{eq:modetranslation}, the same is true for phase $\phi$, and
\begin{equation}\label{eq:phase}
\widehat{\Psi_{n,\phi}^\nu e_{k,\ell}}(k,\ell+pN)
=e^{ip\phi}\widehat{\Psi_{n,0}^\nu e_{k,\ell}}(k,\ell+pN).
\end{equation}

Consider a cross term in \eqref{eq:crossterms} with $\ell\neq\ell'$. Both outgoing coefficients can be nonzero only if
\[
q=\ell+pN=\ell'+p'N
\qquad\text{for some }p,p'\in\Z.
\]
Since $\ell\neq\ell'$, we have $p\neq p'$. By \eqref{eq:phase}, this cross term equals its value at $\phi=0$ multiplied by $e^{i(p-p')\phi}$. Hence its phase average vanishes:
\[
\frac1{2\pi}\int_0^{2\pi}e^{i(p-p')\phi}\,\dee\phi=0.
\]
Summing over the finitely many input pairs and all output frequencies in \eqref{eq:crossterms} proves the lemma.
\end{proof}

\begin{proof}[Proof of Lemma~\ref{lem:transfer}]
Averaging \eqref{eq:leakageexpansion} and applying Lemma~\ref{lem:crossaverage} gives
\begin{equation}\label{eq:summedaverage}
\frac1{2\pi}\int_0^{2\pi}
\bigl\|(1-P_{n+1})\Psi_{n,\phi}^\nu g\bigr\|_{L^2}^2\,\dee\phi
=\frac1{2\pi}\int_0^{2\pi}\mathcal D_\phi\,\dee\phi
\leqc\frac1m\left(\frac4K\right)^n\|g\|_{L^2}^2,
\end{equation}
where the inequality follows from \eqref{eq:diagonalbound}. This proves \eqref{eq:transfer} for even $n$. The case of odd $n$ follows by interchanging the coordinates.
\end{proof}

\subsection{Selection of the phase}\label{sec:crossvariation}

For each fixed $\nu\in[0,1]$, the phase-averaged estimate
\eqref{eq:transfer} directly yields a phase $\phi=\phi(\nu)$ such that
\[
\bigl\|(1-P_{n+1})\Psi_{n,\phi(\nu)}^\nu g_n^\nu\bigr\|_{L^2}^2
\le\frac C m\left(\frac4K\right)^n\|g_n^\nu\|_{L^2}^2
\le\frac C m\left(\frac4K\right)^n\|f_0\|_{L^2}^2.
\]
However, this phase may depend on $\nu$, both through the solution map and through the incoming scalar $g_n^\nu$. The goal of this subsection is to choose a single phase $\phi_n$, independent of $\nu$, for which the uniform bound \eqref{eq:phasechoice} holds, thereby proving Lemma~\ref{lem:phasechoice}.

The phases are chosen inductively. At step $n$, the phases $\phi_0,\ldots,\phi_{n-1}$ selected at the preceding steps, independently of $\nu$, determine the incoming family $\{g_n^\nu\}_{\nu\in[0,1]}$ through \eqref{eq:iteration}, starting from $g_0^\nu=f_0$. We regard this family as given while analyzing step $n$. For any phase $\phi$, denote the Fourier leakage norm by
\[
\mathrm E_{n,\phi}^\nu = \bigl\|(1-P_{n+1})\Psi_{n,\phi}^\nu g_n^\nu\bigr\|_{L^2}.
\]
By Lemma~\ref{lem:transfer} and $\|g_n^\nu\|_{L^2}\le\|f_0\|_{L^2}$, the phase average of its square satisfies
\begin{equation}\label{eq:leakageaverage}
\frac1{2\pi}\int_0^{2\pi}\bigl(\mathrm E_{n,\phi}^\nu\bigr)^2\,\dee\phi
\le\frac C m\left(\frac4K\right)^n\|f_0\|_{L^2}^2.
\end{equation}

To choose a phase that works for all $\nu \in [0,1]$, we study how $\mathrm E_{n,\phi}^\nu$ varies with $\nu$ while keeping $\phi$ fixed. The idea is that small leakage at one diffusivity persists at other values provided the accumulated variation is also small. The following lemma gives quantitative control of this variation. Combined with the phase-averaged bound \eqref{eq:leakageaverage}, it will allow us to control the phase average of $\sup_{0\le\nu\le1}(\mathrm E_{n,\phi}^\nu)^3$ and hence choose a single phase with uniformly small leakage in the proof of Lemma~\ref{lem:phasechoice}.

\begin{lemma}\label{lem:leakagevariation}
For every phase $\phi$, the map $\nu\mapsto\mathrm E_{n,\phi}^\nu$ is locally Lipschitz on $(0,1]$ and satisfies, for almost every $\nu\in(0,1)$,
\begin{equation}\label{eq:leakagevariation}
\left|\partial_\nu\mathrm E_{n,\phi}^\nu\right|
\le w_n(\nu)\|f_0\|_{L^2},\qquad\text{with weight}\qquad 
w_n(\nu) = \min\left\{\frac1\nu,\frac{me^{K(n+1)}}{\sqrt\nu}\right\}.
\end{equation}
\end{lemma}

\begin{proof}
We give an a priori calculation by formally differentiating in $\nu$. This can be justified by a standard difference-quotient argument.

Let $g^\nu$ be the evolution obtained by solving \eqref{eq:ADE} between successive projection times, starting from $g^\nu(0)=f_0$. On $[t_j,t_{j+1})$, its initial datum is $g^\nu(t_j)=g_j^\nu$; the shear has phase $\phi_j$ for $j<n$ and phase $\phi$ for $j=n$. We denote the terminal value before projection by $g^\nu(t_{j+1}^-)$. By \eqref{eq:iteration},
\[
g^\nu(t_{j+1})=g_{j+1}^\nu=P_{j+1}g^\nu(t_{j+1}^-),
\qquad 0\le j<n.
\]
Thus $g^\nu$ may jump at the intermediate projection times. Since each $P_{j+1}$ is a Fourier projection, it commutes with derivatives and is a contraction in both $L^2$ and $\dot H^1$. Consequently, the $L^2$ norms of $g^\nu$, $\partial_\nu g^\nu$, and $\grad g^\nu$ do not increase at the intermediate projection times $t_{j+1}$.

At the final time, the leakage is
\[
\mathrm E_{n,\phi}^\nu
=\bigl\|(1-P_{n+1})g^\nu(t_{n+1}^-)\bigr\|_{L^2}.
\]
The reverse triangle inequality and the $L^2$ contraction of $1-P_{n+1}$ give, for almost every $\nu\in(0,1)$,
\begin{equation}\label{eq:leakagevariationsum}
\left|\partial_\nu\mathrm E_{n,\phi}^\nu\right|
\le\bigl\|(1-P_{n+1})\partial_\nu g^\nu(t_{n+1}^-)\bigr\|_{L^2}
\le\|\partial_\nu g^\nu(t_{n+1}^-)\|_{L^2}.
\end{equation}

It remains to control $\|\partial_\nu g^\nu(t_{n+1}^-)\|_{L^2}$. Since the shears and projections are independent of the diffusivity, $\partial_\nu g^\nu$ satisfies
\begin{equation}\label{eq:viscosityderivativepde}
\partial_t(\partial_\nu g^\nu)+u_j\cdot\nabla(\partial_\nu g^\nu)
=\nu\Delta(\partial_\nu g^\nu)+\Delta g^\nu,
\qquad \partial_\nu g^\nu(t_j)=\partial_\nu g_j^\nu,
\end{equation}
on $[t_j,t_{j+1})$, with
\[
\partial_\nu g_0^\nu=0,\qquad\text{and}\qquad
\partial_\nu g_{j+1}^\nu
=P_{j+1}\partial_\nu g^\nu(t_{j+1}^-),\qquad 0\leq j<n.
\]
Testing \eqref{eq:viscosityderivativepde} against $\partial_\nu g^\nu$ and applying Young's inequality yields, on each interval,
\begin{equation}\label{eq:viscosityderivativeenergy}
\frac{\dee}{\dee t}\|\partial_\nu g^\nu\|_{L^2}^2
+\nu\|\grad\partial_\nu g^\nu\|_{L^2}^2
\le\frac1\nu\|\grad g^\nu\|_{L^2}^2.
\end{equation}

Recall the energy identity
\[
\frac{\dee}{\dee t}\|g^\nu\|_{L^2}^2
+2\nu\|\grad g^\nu\|_{L^2}^2=0.
\]
This implies
\[
\frac{\dee}{\dee t}\left(
\|\partial_\nu g^\nu\|_{L^2}^2
+\frac1{2\nu^2}\|g^\nu\|_{L^2}^2
\right)\le0
\]
on each interval. Together with contraction at the intermediate projection times, we obtain
\begin{align*}
\|\partial_\nu g^\nu(t_{n+1}^-)\|_{L^2}^2
&\le\|\partial_\nu g^\nu(t_{n+1}^-)\|_{L^2}^2 + \frac1{2\nu^2}\|g^\nu(t_{n+1}^-)\|_{L^2}^2\\
&\le\|\partial_\nu g^\nu(0)\|_{L^2}^2 +\frac1{2\nu^2}\|g^\nu(0)\|_{L^2}^2
=\frac1{2\nu^2}\|f_0\|_{L^2}^2.
\end{align*}
Substituting into \eqref{eq:leakagevariationsum} gives the bound in \eqref{eq:leakagevariation}, with weight $1/\nu$.
Since this weight is not integrable at $\nu=0$, we obtain another bound by studying the dynamics of $\grad g^\nu$. The gradient inequality
\[
\frac{\dee}{\dee t}\|\grad g^\nu\|_{L^2}^2
\le2\alpha_jN_j\|\grad g^\nu\|_{L^2}^2,
\]
together with contraction at projection times and $\alpha_jN_j\tau_j=K$, gives by Gronwall's inequality
\begin{equation}\label{eq:inputgradient}
\|\grad g^\nu(t)\|_{L^2}
\le me^{K(n+1)}\|f_0\|_{L^2},
\qquad 0\le t<t_{n+1}.
\end{equation}
Returning to \eqref{eq:viscosityderivativeenergy}, we obtain
\[
\frac{\dee}{\dee t}\|\partial_\nu g^\nu\|_{L^2}^2
\le\frac{m^2e^{2K(n+1)}}\nu\|f_0\|_{L^2}^2.
\]
Since $\partial_\nu g^\nu$ starts from zero, its $L^2$ norm does not increase at projection times, and $t_{n+1}\le1$, this yields
\[
\|\partial_\nu g^\nu(t_{n+1}^-)\|_{L^2}
\le\frac{me^{K(n+1)}}{\sqrt\nu}\|f_0\|_{L^2}.
\]
Combining the two bounds in \eqref{eq:leakagevariationsum} proves \eqref{eq:leakagevariation}. These bounds and the reverse triangle inequality also give the asserted local Lipschitz continuity.
\end{proof}

Note that the weight $w_n(\nu)$ is integrable on $(0,1)$. Indeed, we have the explicit calculation:
\begin{equation}\label{eq:variationweight}
\int_0^1 w_n(\nu)\,\dee\nu
=\int_0^{m^{-2}e^{-2K(n+1)}}\frac{me^{K(n+1)}}{\sqrt\nu}\,\dee\nu
+\int_{m^{-2}e^{-2K(n+1)}}^1\frac1\nu\,\dee\nu
=2+2\log m+2K(n+1).
\end{equation}

We are now ready to prove Lemma~\ref{lem:phasechoice}.

\begin{proof}[Proof of Lemma~\ref{lem:phasechoice}]
By \eqref{eq:leakagevariation} and the chain rule,
\[
\left|\partial_\nu\bigl(\mathrm E_{n,\phi}^\nu\bigr)^3\right|
\le3\|f_0\|_{L^2}w_n(\nu)\bigl(\mathrm E_{n,\phi}^\nu\bigr)^2.
\]
Integrating over $[\nu_0,1]$ gives, for every phase $\phi$ and every $0<\nu_0\le1$,
\[
\bigl(\mathrm E_{n,\phi}^{\nu_0}\bigr)^3
\le\bigl(\mathrm E_{n,\phi}^1\bigr)^3
+3\|f_0\|_{L^2}\int_{\nu_0}^1w_n(\nu)\bigl(\mathrm E_{n,\phi}^\nu\bigr)^2\,\dee\nu.
\]
Since $\mathrm E_{n,\phi}^\nu\le\|f_0\|_{L^2}$ and $w_n$ is integrable by \eqref{eq:variationweight}, we may replace the lower endpoint by zero and, by continuity at $\nu=0$, take the supremum over $[0,1]$ to obtain
\[
\sup_{\nu\in[0,1]}\bigl(\mathrm E_{n,\phi}^\nu\bigr)^3
\le\bigl(\mathrm E_{n,\phi}^1\bigr)^3
+3\|f_0\|_{L^2}\int_0^1w_n(\nu)\bigl(\mathrm E_{n,\phi}^\nu\bigr)^2\,\dee\nu.
\]
Averaging in $\phi$, using $(\mathrm E_{n,\phi}^1)^3\le\|f_0\|_{L^2}(\mathrm E_{n,\phi}^1)^2$, and applying \eqref{eq:leakageaverage} and \eqref{eq:variationweight}, we obtain
\begin{align*}
\frac1{2\pi}\int_0^{2\pi}\sup_{\nu\in[0,1]}
\bigl(\mathrm E_{n,\phi}^\nu\bigr)^3\,\dee\phi
&\le\frac{\|f_0\|_{L^2}}{2\pi}\int_0^{2\pi}\bigl(\mathrm E_{n,\phi}^1\bigr)^2\,\dee\phi
+3\|f_0\|_{L^2}\int_0^1w_n(\nu)\left(\frac1{2\pi}\int_0^{2\pi}\bigl(\mathrm E_{n,\phi}^\nu\bigr)^2\,\dee\phi\right)\dee\nu\\
&\le\frac C m\left(\frac4K\right)^n\|f_0\|_{L^2}^3\left(1+\int_0^1w_n(\nu)\,\dee\nu\right)\\
&\le\frac C m\left(\frac4K\right)^n\bigl(1+\log m+K(n+1)\bigr)\|f_0\|_{L^2}^3.
\end{align*}
Choose a phase $\phi_n$ for which the supremum does not exceed its phase average. Taking cube roots gives
\[
\sup_{\nu\in[0,1]} \bigl\|(1-P_{n+1})\Psi_{n,\phi_n}^\nu g_n^\nu\bigr\|_{L^2}
\le C\left(\frac{1+\log m+K(n+1)}m\right)^{1/3} \left(\frac4K\right)^{n/3}\|f_0\|_{L^2}.
\]
Choosing $K$ large enough, the coefficient decays geometrically in $n$. In particular, for $K>32$ we have the explicit bound
\[
C\left(\frac{1+\log m+K(n+1)}m\right)^{1/3}\left(\frac4K\right)^{n/3}
\le C_K\left(\frac{1+\log m}{m}\right)^{1/3}\,2^{-n}.
\]
The constant $C_K$ depends only on $K$ and is independent of $n$ and the previous phases. Finally, choose one $m$ large enough, depending on $K$, such that $C_K((1+\log m)/m)^{1/3}\le1/4$. This proves \eqref{eq:phasechoice} for $0\le\nu\le1$ at every step.
\end{proof}

\section{Proof of Theorem~\ref{thm:main}}\label{sec:concludeproof}

We now prove Theorem~\ref{thm:main} by combining the Fourier leakage bound with dissipation at high frequencies. Lemma~\ref{lem:phasechoice} provides a single sequence of phases, independent of viscosity, for which the errors introduced by the successive Fourier projections are summable. Consequently, the projected scalars $g_n^\nu$ approximate the actual solution $f^\nu(t_n)$ with an error bounded uniformly in both $n$ and $\nu$.

As in the Couette evolution, the Fourier support of $g_n^\nu$ lies in $C_n$ and therefore moves to arbitrarily high frequencies as $t_n\to1$. We show that $\|g_n^\nu\|_{L^2}\to0$ for each fixed $\nu>0$. Together with the uniform approximation bound, this yields a positive loss of energy for the actual solution that persists as $\nu\to0$, establishing anomalous dissipation.

\begin{proof}[Proof of Theorem~\ref{thm:main}]

Fix $K>32$ and choose $m\in\N$ and the $\nu$-independent phases $\{\phi_n\}_{n\ge0}$ provided by Lemma~\ref{lem:phasechoice}, so that \eqref{eq:phasechoice} holds for every $n\ge0$. This results in a $\nu$-independent velocity field $u$ with the regularity established in Section~\ref{sec:osgood}. Since the phases are now fixed, to simplify notation we write $\Psi_{n}^\nu$ instead of $\Psi_{n,\phi_n}^\nu$ in the remainder of this proof.

Fix $0<\nu\le1$. Lemma~\ref{lem:phasechoice} bounds the error introduced by the Fourier projection at each step. Summing these errors gives a bound on the total approximation error between $g_n^\nu$ and $f^\nu(t_n)$, uniform in both $n$ and $\nu$. Let 
\[
E_n = \|f^\nu(t_n) - g_n^\nu\|_{L^2}.
\]
By the triangle inequality, the choice of phases satisfying \eqref{eq:phasechoice}, and $\|\Psi_{n}^\nu\|_{L^2\to L^2} \le 1$, we have 
\begin{align*}
    E_{n+1} & = \|\Psi_n^\nu f^\nu(t_n) - P_{n+1}\Psi_n^\nu g_n^\nu\|_{L^2} \\ 
    & \le \|\Psi_n^\nu(f^\nu(t_n) - g_n^\nu)\|_{L^2} + \|(1-P_{n+1})\Psi_n^\nu g_n^\nu\|_{L^2} 
    \le E_n + 2^{-(n+2)}\|f_0\|_{L^2}.
\end{align*}
Since $E_0 = 0$, iterating the above gives
\begin{equation}\label{eq:iterationerror}
    E_n \le \|f_0\|_{L^2} \sum_{k=0}^{n-1}2^{-(k+2)} \le \frac{1}{2}\|f_0\|_{L^2}.
\end{equation}
By \eqref{eq:iterationerror}, for every $n \ge 0$ there holds 
$$ \|f^\nu(t_n)\|_{L^2}^2 \le (E_n + \|g_n^\nu\|_{L^2})^2 \le \frac{1}{2}\|f_0\|_{L^2}^2 + 2\|g_n^\nu\|_{L^2}^2. $$
It then follows from the energy identity that
\begin{equation}\label{eq:dissipationfromapproximation}
2\nu\int_0^1\|\grad f^\nu(t)\|_{L^2}^2\,\dee t
\ge2\nu\int_0^{t_n}\|\grad f^\nu(t)\|_{L^2}^2\,\dee t
=\|f_0\|_{L^2}^2-\|f^\nu(t_n)\|_{L^2}^2
\ge\frac12\|f_0\|_{L^2}^2-2\|g_n^\nu\|_{L^2}^2.
\end{equation}
Therefore, to complete the proof of Theorem~\ref{thm:main}, it suffices to show that 
\begin{equation}\label{eq:approxdecay}
    \lim_{n \to \infty} \|g_n^\nu\|_{L^2} = 0
\end{equation}
for each fixed $0<\nu\le1$. Indeed, assuming \eqref{eq:approxdecay} and letting $n\to\infty$ in \eqref{eq:dissipationfromapproximation}, we obtain
\[
\nu\int_0^1\|\grad f^\nu(t)\|_{L^2}^2\,\dee t
\ge\frac14\|f_0\|_{L^2}^2,
\qquad 0<\nu\le1.
\]
The lower bound is positive and independent of $\nu$, so taking the lower limit as $\nu\to0$ gives the conclusion of Theorem~\ref{thm:main}.

It remains to show \eqref{eq:approxdecay}, that is, complete dissipation of the auxiliary scalar for each fixed $\nu>0$. Recall the evolution $g^\nu$ defined in the proof of Lemma~\ref{lem:leakagevariation}, with the phases now fixed as above. Fix $n\ge0$ and consider its restriction to $[t_n,t_{n+1})$, starting from $g_n^\nu$. Suppose without loss of generality that $n$ is even. Then, $g_n^\nu$ has Fourier support in $C_n$ with
$$C_n\subset\bigl\{(k,\ell)\in\Z^2: |k|\ge m(K/2)^n\}.$$
Since advection-diffusion with a horizontal shear flow preserves $x$-frequencies, we have
$$
\|\grad g^\nu(t)\|_{L^2}^2 \ge \|\partial_x g^\nu(t)\|_{L^2}^2
\ge m^2(K/2)^{2n}\|g^\nu(t)\|_{L^2}^2
$$
for all $t \in [t_n,t_{n+1})$. Therefore, 
$$ \frac{\dee}{\dee t} \|g^\nu(t)\|_{L^2}^2 = -2\nu\|\grad g^\nu(t)\|_{L^2}^2 \le -2\nu m^2 (K/2)^{2n}\|g^\nu(t)\|_{L^2}^2. $$
Applying Gr\"{o}nwall's inequality up to $t_{n+1}^-$ and then using the contraction of $P_{n+1}$ yields
\begin{equation}\label{eq:onestepapproxdecay}
\|g^\nu_{n+1}\|_{L^2} \le \|g^\nu(t_{n+1}^-)\|_{L^2} \le e^{-\nu \tau_n m^2 (K/2)^{2n}}\|g_n^\nu\|_{L^2}.
\end{equation}
The same argument shows that \eqref{eq:onestepapproxdecay} also holds for odd $n$. Thus, we may iterate the bound to obtain 
\begin{equation}\label{eq:gngronwall}
\|g_n^\nu\|_{L^2}
\le
\|f_0\|_{L^2}\exp\!\left(
-\nu m^2\sum_{j=0}^{n-1}\tau_j\left(\frac{K}{2}\right)^{2j} 
\right).
\end{equation}
With the choice of $\tau_n$ in \eqref{eq:taundef}, we have 
\begin{equation}\label{eq:diverge}
\sum_{j=0}^\infty\tau_j\left(\frac{K}{2}\right)^{2j} = 
\sum_{j=0}^\infty\frac Kc\frac{2^{-j}}{\omega(\pi 2^{-j})}\left(\frac{K}{2}\right)^{2j} 
\ge \frac{K}{c\omega(\pi)}
\sum_{j=0}^\infty\left(\frac{K^2}{8}\right)^j
=\infty.
\end{equation}
Together with \eqref{eq:gngronwall}, this proves \eqref{eq:approxdecay} for each fixed $0<\nu\le1$.

\end{proof}

\bigskip
\subsection*{Disclosure of AI assistance}
The authors used ChatGPT 6 Astra (OpenAI) to assist in developing some of the key ideas and proof strategies in this work, as well as in revising the manuscript. The authors independently verified the mathematical arguments and take full responsibility for the results and the final manuscript.

\bigskip

\bigskip

\addcontentsline{toc}{section}{References}

\begin{thebibliography}{10}

\bibitem{ACM}
Giovanni Alberti, Gianluca Crippa, and Anna~L. Mazzucato.
\newblock Exponential self-similar mixing by incompressible flows.
\newblock {\em Journal of the American Mathematical Society}, 32(2):445--490, 2019.

\bibitem{Ambrosio2004}
Luigi Ambrosio.
\newblock Transport equation and {Cauchy} problem for {BV} vector fields.
\newblock {\em Inventiones Mathematicae}, 158(2):227--260, 2004.

\bibitem{AmbrosioBernard2008}
Luigi Ambrosio and Patrick Bernard.
\newblock Uniqueness of signed measures solving the continuity equation for
  {Osgood} vector fields.
\newblock {\em Rendiconti Lincei. Matematica e Applicazioni}, 19(3):237--245, 2008.

\bibitem{AmbrosioCrippa2014}
Luigi Ambrosio and Gianluca Crippa.
\newblock Continuity equations and {ODE} flows with non-smooth velocity.
\newblock {\em Proceedings of the Royal Society of Edinburgh Section A: Mathematics}, 144(6):1191--1244, 2014.

\bibitem{AV}
Scott Armstrong and Vlad Vicol.
\newblock Anomalous diffusion by fractal homogenization.
\newblock {\em Annals of PDE}, 11(1):2, 2025.

\bibitem{BonicattoCiampaCrippa2022}
Paolo Bonicatto, Gennaro Ciampa, and Gianluca Crippa.
\newblock On the advection-diffusion equation with rough coefficients: weak
  solutions and vanishing viscosity.
\newblock {\em Journal de Math{\'e}matiques Pures et Appliqu{\'e}es},
  167:204--224, 2022.

\bibitem{BonicattoGusev2019}
Paolo Bonicatto and Nikolay~A. Gusev.
\newblock Non-uniqueness of signed measure-valued solutions to the continuity
  equation in presence of a unique flow.
\newblock {\em Rendiconti Lincei. Matematica e Applicazioni}, 30(3):511--531, 2019.

\bibitem{BCCDS}
Elia Bru{\`e}, Maria Colombo, Gianluca Crippa, Camillo De~Lellis, and Massimo Sorella.
\newblock Onsager critical solutions of the forced {Navier--Stokes} equations.
\newblock {\em Communications on Pure and Applied Analysis}, 23(10):1350--1366,
  2024.

\bibitem{BDL}
Elia Bru{\`e} and Camillo De~Lellis.
\newblock Anomalous dissipation for the forced {3D} {Navier--Stokes} equations.
\newblock {\em Communications in Mathematical Physics}, 400(3):1507--1533,
  2023.

\bibitem{BurczakSzekelyhidiWu2023}
Jan Burczak, L{\'a}szl{\'o} Sz{\'e}kelyhidi, Jr., and Bian Wu.
\newblock Anomalous dissipation and {Euler} flows.
\newblock arXiv preprint arXiv:2310.02934, 2023.

\bibitem{CCS}
Maria Colombo, Gianluca Crippa, and Massimo Sorella.
\newblock Anomalous dissipation and lack of selection in the {Obukhov--Corrsin}
  theory of scalar turbulence.
\newblock {\em Annals of PDE}, 9(2):21, 2023.

\bibitem{ColomboKumar2026}
Roberto Colombo and Anuj Kumar.
\newblock Sharpness of the {Osgood} criterion for the continuity equation with
  divergence-free vector fields.
\newblock arXiv preprint arXiv:2601.12096, 2026.

\bibitem{Corrsin51}
Stanley Corrsin.
\newblock On the spectrum of isotropic temperature fluctuations in an isotropic
  turbulence.
\newblock {\em Journal of Applied Physics}, 22(4):469--473, 1951.

\bibitem{Depauw2003}
Nicolas Depauw.
\newblock Non unicit{\'e} des solutions born{\'e}es pour un champ de vecteurs
  {BV} en dehors d'un hyperplan.
\newblock {\em Comptes Rendus Math{\'e}matique}, 337(4):249--252, 2003.

\bibitem{DL}
Ronald~J. DiPerna and Pierre-Louis Lions.
\newblock Ordinary differential equations, transport theory and {Sobolev}
  spaces.
\newblock {\em Inventiones Mathematicae}, 98(3):511--547, 1989.

\bibitem{DSrY05}
Diego~A. Donzis, Katepalli~R. Sreenivasan, and Pui-Kuen Yeung.
\newblock Scalar dissipation rate and dissipative anomaly in isotropic
  turbulence.
\newblock {\em Journal of Fluid Mechanics}, 532:199--216, 2005.

\bibitem{DEIJ22}
Theodore~D. Drivas, Tarek~M. Elgindi, Gautam Iyer, and In-Jee Jeong.
\newblock Anomalous dissipation in passive scalar transport.
\newblock {\em Archive for Rational Mechanics and Analysis}, 243(3):1151--1180, 2022.

\bibitem{ElgindiLiss2024}
Tarek~M. Elgindi and Kyle~L. Liss.
\newblock Norm growth, non-uniqueness, and anomalous dissipation in passive
  scalars.
\newblock {\em Archive for Rational Mechanics and Analysis}, 248(6):120, 2024.

\bibitem{EyinkSreenivasan2006}
Gregory~L. Eyink and Katepalli~R. Sreenivasan.
\newblock {Onsager} and the theory of hydrodynamic turbulence.
\newblock {\em Reviews of Modern Physics}, 78(1):87--135, 2006.

\bibitem{Frisch}
Uriel Frisch.
\newblock {\em Turbulence: The Legacy of {A. N. Kolmogorov}}.
\newblock Cambridge University Press, 1995.

\bibitem{HessChildsRowan2025}
Elias Hess-Childs and Keefer Rowan.
\newblock A universal total anomalous dissipator.
\newblock {\em Communications in Mathematical Physics}, 407:257, 2026.

\bibitem{JS2d}
Carl Johan Peter Johansson and Massimo Sorella.
\newblock Anomalous dissipation via spontaneous stochasticity with a
  two-dimensional autonomous velocity field.
\newblock {\em Duke Mathematical Journal}, 175(11):1957--2024, 2026.

\bibitem{JS}
Carl Johan Peter Johansson and Massimo Sorella.
\newblock Nontrivial absolutely continuous part of anomalous dissipation
  measures in time.
\newblock {\em Journal of Differential Equations}, 453(4):113912, 2026.

\bibitem{K41a}
Andrey Nikolaevich Kolmogorov.
\newblock Dissipation of energy in the locally isotropic turbulence.
\newblock {\em Doklady Akademii Nauk SSSR}, 32:19--21, 1941.

\bibitem{Obukhov49}
Alexander Mikhailovich Obukhov.
\newblock Structure of the temperature field in a turbulent flow.
\newblock {\em Izvestiya Akademii Nauk SSSR. Seriya Geograficheskaya i Geofizicheskaya}, 13(1):58--69, 1949.

\bibitem{Osgood1898}
William Fogg Osgood.
\newblock Beweis der {Existenz} einer {L{\"o}sung} der {Differentialgleichung}
  $dy/dx=f(x,y)$ ohne {Hinzunahme} der {Cauchy-Lipschitz}'schen {Bedingung}.
\newblock {\em Monatshefte f{\"u}r Mathematik und Physik}, 9(1):331--345, 1898.

\bibitem{Rowan2024Kraichnan}
Keefer Rowan.
\newblock On anomalous diffusion in the {Kraichnan} model and
  correlated-in-time variants.
\newblock {\em Archive for Rational Mechanics and Analysis}, 248(5):93, 2024.

\bibitem{SS2000}
Boris~I. Shraiman and Eric~D. Siggia.
\newblock Scalar turbulence.
\newblock {\em Nature}, 405(6787):639--646, 2000.

\bibitem{YaoZlatos2017}
Yao Yao and Andrej Zlato{\v{s}}.
\newblock Mixing and un-mixing by incompressible flows.
\newblock {\em Journal of the European Mathematical Society}, 19(7):1911--1948,
  2017.

\end{thebibliography}
\end{document}